\documentstyle[amssymb,amsfonts]{amsart}

\def\bi{\begin{itemize}}
\def\bs{\begin{split}}
\def\es{\end{split}}
\def\ba{\begin{align}}
\def\bas{\begin{align*}}
\def\ea{\end{align}}
\def\eas{\end{align*}}
\def\R{{\hbox{\bf R}}}

\def\Im{{\hbox{Im}}}
\def\Re{{\hbox{Re}}}

\def\wb{{\hbox{wb}}}
\def\bound{{\hbox{b}}}
\def\R{{\hbox{\bf R}}}
\def\C{{\hbox{\bf C}}}
\def\sgn{{\hbox{sgn}}}

\def\Z{{\hbox{\bf Z}}}

\def\eps{\varepsilon}

\newenvironment{proof}{\noindent {\bf Proof} }{\endprf\par}
\def \endprf{\hfill  {\vrule height6pt width6pt depth0pt}\medskip}
\def\emph#1{{\it #1}}
\def\textbf#1{{\bf #1}}

\theoremstyle{plain}
  \newtheorem{theorem}[subsection]{Theorem}
  
  \newtheorem{proposition}[subsection]{Proposition}
  \newtheorem{lemma}[subsection]{Lemma}
  \newtheorem{corollary}[subsection]{Corollary}

\theoremstyle{remark}

\theoremstyle{definition}

\include{psfig}

\begin{document}

\title[Asymptotic behavior of focusing NLS]{On the asymptotic behavior of large radial data for a focusing non-linear Schr\"odinger equation}
\author{Terence Tao}
\address{Department of Mathematics, UCLA, Los Angeles CA 90095-1555}
\email{tao@@math.ucla.edu}
\thanks{This work was conducted at Australian National University.  The author is a Clay Prize Fellow and is supported by a grant from the Packard foundation.  We thank Jim Colliander and Igor Rodnianski for helpful comments and corrections, and Wilhelm Schlag for pointing out that the author's original proof of Theorem \ref{easy} was incorrect.}
\subjclass{35Q55}

\vspace{-0.3in}
\begin{abstract}
We study the asymptotic behavior of large data radial solutions to the focusing Schr\"odinger equation
$i u_t + \Delta u = -|u|^2 u$ in $\R^3$,
assuming globally bounded $H^1(\R^3)$ norm (i.e. no blowup in the energy space).  We show that as
$t \to \pm \infty$, these solutions split into the sum of three terms: a radiation term that evolves according to the linear Schr\"odinger equation, a smooth function localized near the origin, and an error that goes to zero in the $\dot H^1(\R^3)$ norm.  Furthermore, the smooth function near the origin is either zero (in which case one has scattering to a free solution), or has mass and energy bounded strictly away from zero, and obeys an asymptotic Pohozaev identity.
These results are consistent with the conjecture of soliton resolution.
\end{abstract}

\maketitle

\section{Introduction}

In this paper we consider the (forward-in-time) asymptotic behavior of global energy-class solutions $u: [0,+\infty) \times \R^3 \to \C$ to the focusing cubic\footnote{It is likely that at least some of the results here also extend to other dimensions and exponents.  However, we need to be in three and higher dimensions in order for the fundamental solution $e^{it\Delta}$ to be ``short range'' (decaying faster than $t^{-1}$), and the algebraic nature of the cubic non-linearity is also somewhat convenient.  Also much of the asymptotic analysis is restricted to the $L^2$-supercritical regime $p > 1 + \frac{4}{d}$ in order for global Strichartz estimates to be useful.} non-linear Schr\"odinger equation
\begin{equation}\label{nls}
i u_t + \Delta u = F(u)
\end{equation}
in three dimensions, where $F(u) := - |u|^2 u$ is the focusing cubic non-linearity.

Non-linear Schr\"odinger equations such as \eqref{nls} have been intensively studied and we review\footnote{As this is a vast field, we cannot hope to come close to an exhaustive description of results here, and the references given are not intended to be complete.  We refer the reader to the books \cite{caz}, \cite{sulem}, \cite{borg:book}, \cite{strauss} for a more detailed survey.} some of the results here. The equation \eqref{nls} is locally well-posed in the energy class $H^1(\R^3)$ (see e.g. \cite{caz}, \cite{sulem}, \cite{strauss}; in fact it is locally well-posed in $H^s(\R^3)$ for all $s \geq 1/2$ \cite{cwI}).  It also enjoys the conservation of mass
$$ MASS(u(t)) := \int_{\R^3} |u(t,x)|^2\ dx$$
and conservation of the Hamiltonian
$$ ENERGY(u(t)) := \int_{\R^3} \frac{1}{2} |\nabla u(t,x)|^2 - \frac{1}{4} |u(t,x)|^4\ dx,$$
however these two conservation laws are unfortunately not enough to control the $H^1(\R^3)$ norm; the Gagliardo
Nirenberg inequality 
\begin{equation}\label{gn}
 \int |u|^4 \leq C(\int |\nabla u|^2)^{3/2} (\int |u|^2)^{1/2}
\end{equation}
is not sufficient unless the mass\footnote{More precisely, one needs a scale invariant quantity such as $(\int |u|^2)(\int |\nabla u|^2)$ to be small at time zero, and then it will be small for all time by combining the above argument with the continuity method.} is sufficiently small.  Indeed, if the Hamiltonian is negative then one has blowup in finite time from the virial identity \cite{glassey}, \cite{ozawa} (see \cite{caz}, \cite{sulem} for further discussion; see also Corollary \ref{energy} below).

In this paper we will consider the asymptotic behavior (as\footnote{The behavior as $t \to -\infty$ is completely identical, thanks to the time reversal symmetry $u(t,x) \mapsto \overline{u(-t,x)}$ of \eqref{nls}, and will not be discussed here.  It is however an interesting question as to whether the solutions which scatter to a free solution as $t \to +\infty$ match at all with the solutions which scatter to a free solution as $t \to -\infty$; this is true in the completely integrable one-dimensional situation, as well as for small data, but may well be false for the three-dimensional large data case.} $t \to +\infty$) of solutions which do not blow up in the energy class; more precisely
we consider solutions where we assume \emph{a priori} the uniform $H^1$ bound
\begin{equation}\label{energy-bound}
\sup_{t \in [0,+\infty)} \| u(t) \|_{H^1(\R^3)} \leq E
\end{equation}
for some $0 < E < +\infty$ (note we do not assume any smallness assumption on $E$; we also allow the possibility that
solution could blow up at some finite \emph{negative} time).  This condition has appeared in other work on this equation, see in particular \cite{borg:growth}.  In particular (by \eqref{gn}) we assume that the mass and energy of $u$ is finite.

Bounded energy solutions will occur for instance when the product of the mass and energy of $u$ (which is a scale-invariant quantity) is sufficiently small, thanks to the Gagliardo-Nirenberg inequality \eqref{gn}; another
option is to modify the non-linearity $F(u)$ for large $u$ by replacing it smoothly with a non-linearity which
behaves like $O(|u|^p)$ for some $p < 1 + \frac{4}{3}$, as this allows the Gagliardi-Nirenberg inequality argument
to derive \eqref{energy-bound} even when the energy and mass are large (as we shall see, it is the low values of $u$ which will dominate our discussion, the high values being irrelevant except in order to establish a local existence theory).

With the \emph{a priori} assumption \eqref{energy-bound} we know that $u$ is globally well-posed in $H^1(\R^3)$ by iterating
the local well-posedness argument.  But this standard global existence result does not reveal very much about the asymptotic behavior of
$u$.  For data whose energy $E$ is small it is possible to show that the solution $u(t)$ 
eventually scatters to approach a free solution $e^{it\Delta} u_+$ in the energy norm $H^1(\R^3)$
as $t \to + \infty$, see e.g. \cite{caz}, \cite{strauss} (see also Theorem \ref{easy} below); note that the
potential $|u|^2$ in $F(u)$ is a ``short-range'' potential (it decays integrably in time if $\|u(t)\|_\infty$ has the
expected decay of $t^{-3/2}$) and thus does not cause any non-linear corrections (see \cite{ozawa} for more discussion).  Furthermore, the map $u(0) \mapsto u_+$ is a local homeomorphism in the energy space near the origin (see e.g. \cite{strauss}).

For large data we do not expect this scattering behavior\footnote{One can show however that given any large scattering data $u_+ \in H^1(\R^3)$, one can find a time interval $(T_*,+\infty)$ and an energy class
solution $u: (T_*,+\infty) \times \R^3$ such that $\lim_{t \to +\infty}  \| u(t) - e^{it\Delta} u_+ \|_{H^1(\R^3)} = 0$,
basically by solving the Cauchy problem backwards in time from infinity; see \cite{caz}.  But these scattering solutions do not capture all the large data solutions, for instance the soliton
solutions are clearly not in this class.  On the other hand, in the \emph{defocusing} case every large energy (and even certain infinite energy) data scatters to a free solution, while conversely every free solution is the asymptotic limit of a nonlinear solution; see \cite{gv:scatter}, \cite{borg:scatter}, \cite{ckstt:scatter}.  The arguments in the defocusing case rely on Morawetz inequalities, which have an unfavorable sign in the focusing case; it seems that Morawetz inequalities should still have some value for the focusing case but we were unable to obtain any natural application of them in this setting.} to the focusing equation \eqref{nls}, even assuming the energy bound \eqref{energy-bound}, because of the existence of soliton solutions 
(both ground states and excited states).  These soliton solutions are generated by finite-energy solutions $u(x)$ to the non-linear eigenfunction equation
\begin{equation}\label{nl-eigen}
 -\omega u + \Delta u = F(u)
\end{equation}
for some $\omega > 0$ (which is related to, but not quite the same, quantity as the energy $E$ 
in \eqref{energy-bound}). Such solutions are known to be smooth and exponentially rapidly decreasing (see e.g.
\cite{caz}), and if one makes the further assumption that $u$ is non-negative and spherically symmetric then there is a unique
solution to \eqref{nl-eigen} for each $\omega > 0$ \cite{coff}, \cite{strauss}, \cite{blions}; we refer to this $u$ as the \emph{ground state}.
There also exist radial solutions which change sign, see \cite{bgk}; we refer to these as \emph{excited states}.  Note
that if $u$ is either a ground state or excited state then $u(t,x) := e^{i\omega t} u(x)$ is a solution to \eqref{nls}; we refer to these solutions as the \emph{ground soliton} and \emph{excited soliton} solutions
respectively.  These special solutions can be of course modified by the various symmetries of the Schr\"odinger equation such as
scaling, translations, Gallilean transformations, and phase rotation to produce solitons of various velocities, energies,
widths, etc.

The analysis of solutions close to a soliton state is by now well understood; the ground state soliton for the cubic NLS \eqref{nls} is unstable due to the existence of nearby blowup solutions \cite{bc}, \cite{shatah}, although if one mollifies the growth of the nonlinearity $F(u) = |u|^2 u$ near infinity enough then one can regain orbital stability \cite{wein1}, \cite{wein2}, \cite{caz-lions}, \cite{gss}, \cite{gss-2} and asymptotic stability \cite{bp}, \cite{bp2}, \cite{bs}, \cite{cuccagna}, \cite{cuccagna-2} for the ground state given appropriate conditions on the nonlinearity; see \cite{caz} or \cite{schlag} for further discussion.  For perturbations of excited solitons there seems to be far less known; one possible conjecture is that generic perturbations of the excited soliton should lead either to blowup, or to collapse to a less energetic soliton (or to the vacuum state), plus radiation.  More recently there has been some results on the analysis of perturbations of multisoliton solutions, i.e. superpositions of widely separated and receding solitons (see \cite{schlag}, \cite{perelman}; a similar result for the generalized KdV equation is in \cite{mmt}); while these results do not apply directly to the $L^2$-supercritical cubic equation \eqref{nls}, they do apply to certain mollified versions of the equation (see \cite{schlag} for further discussion).  We should also mention the work in \cite{tsai-yau}, \cite{tsai} on non-linear perturbations of ground and excited states of a Schr\"odinger equation with time-independent potential.  These results however rely either on the variational characterization of the soliton one is perturbing, or on the fact that the time-dependent Hamiltonian $-\Delta - |u|^2$ is well approximated by an explicit linear Hamiltonian (e.g. the charge transfer model, or the linearization of \eqref{nls} around a soliton) for which the spectral structure is well understood.  In the general case we do not have good spectral control on this time-dependent Hamiltonian and so these techniques do not seem to extend to general large data.

However, in the special\footnote{We should caution however that the one-dimensional cubic NLS is $L^2$-subcritical, as opposed to the three-dimensional NLS which is $L^2$-supercritical; this seems to play a decisive role in the analysis.  For $L^2$-subcritical equations one has excellent short-time control (in particular, there is no blowup), but very poor long-time control even in the defocusing case, whereas for $L^2$-supercritical focusing equations the situation is reversed.  The $L^2$ critical equation (such as the 2D cubic NLS) is of course delicate at both short times and long times; see \cite{wein1}, \cite{merle1}, \cite{merle2}, \cite{sulem-break} for some typical results regarding this equation.} case of the cubic focusing NLS equation $iu_t + u_{xx} = -2|u|^2 u$ in one dimension, which is completely integrable, one can use the methods of inverse scattering theory to analyze large data.  For instance, if the initial data is smooth and rapidly decreasing, then it is known (see e.g. \cite{zakharov}, \cite{segur}) that the solution eventually resolves into a  finite number of solitons\footnote{This is an oversimplification.  In the case when two or more solitons have equal speeds, it is possible for ``breather'' solutions which are periodic or even quasiperiodic, and localized for all time to appear.  It is even possible to have two solitons of equal speeds and equal heights to recede from each other at a logarithmic rate $\log t$ in one dimension, although this latter phenomenon may be also due to the slow decay of the fundamental solution in one dimension and we do not know if it occurs in three dimensions.  Nevertheless the generic behavior is that of resolution into receding solitons of different speeds, plus radiation, and even in the exceptional cases the non-radiation part of the solution is still very well localized and smooth, and has asymptotically constant velocity.  We thank Jim Colliander for pointing out this subtlety.} receding from each other, plus some dispersive radiation which goes to zero (although
this radiation does not quite converge to a linear solution because in one dimension the cubic power is on the
borderline between long-range and short-range; see \cite{ozawa}).  

Thus one might conjecture that one has similar behavior (resolution into a multisoliton, plus dispersive radiation,
for generic large data (with perhaps some exceptional channels such as breathers or slowly divergent solitons for a small unstable set of data), or else blowup in $H^1(\R^3)$) for other non-linear Schr\"odinger equations such as the ones studied in this paper; we shall refer to this rather informal statement as the \emph{soliton resolution conjecture}.  Apart from the results mentioned above when the data is close to a soliton or multi-soliton, there are only a few results known to support this conjecture. For instance, some progress has been made on controlling the growth of higher Sobolev norms for such equations, see \cite{borg:growth}, \cite{borg:scatter}, \cite{staff:growth}.  In particular for $H^s$ solutions to \eqref{nls} with the bounded energy assumption \eqref{energy-bound}, it is known that the $H^s$ norms grow at most polynomially in time and that the local $H^s$ norm stays bounded; these results are consistent with the phenomenon of resolution into solitons.  The growth of norms between $L^2$ and $H^1$ in the infinite defocusing energy case have also been studied, see \cite{ckstt:7}, \cite{ckstt:8}, again with polynomial type upper bounds on the growth.  But to obtain control on asymptotic behavior uniformly in time (as opposed to polynomial growth bounds) seems quite difficult.  In particular, the soliton resolution conjecture for non-integrable equations and generic large data seems well out of 
reach of current techniques.  

The case of general large data seems too difficult to analyze at present, so we now restrict our attention to the much simpler \emph{spherically symmetric} case, where $u(0,x) = u(0,|x|)$ (and hence $u(t,x) = u(t,|x|)$, by the rotational symmetry of \eqref{nls} and the uniqueness theory).  In this 
case the radial Sobolev inequality\footnote{This inequality can be easily deduced by applying the one-dimensional Sobolev inequality $\| v \|_{L^\infty(\R)} \lesssim \| v \|_{H^1(\R)}$ to the function $v(|x|) := |x| u(|x|)$ and using Hardy's inequality $\| |u|/|x| \|_{L^2(\R^3)} \lesssim \| \nabla u \|_{L^2(\R^3)}$.}
\begin{equation}\label{radial-sobolev}
\| |x| u \|_{L^\infty(\R^3)} \lesssim \| u \|_{H^1(\R^3)}
\end{equation}
combined with \eqref{energy-bound} tells us that $u(t)$ decays away from the origin\footnote{This explains our previous
remark that the large values of $u$ are not particularly relevant for our purposes; indeed, we see from \eqref{radial-sobolev} that they only make a difference near the spatial origin.} for all times $t$, and 
thus we only expect at most a single stationary spherically symmetric soliton at the origin, and no other soliton 
behavior anywhere else.  

Unfortunately even  in the spherically symmetric case we were still unable to settle the soliton resolution conjecture.
However, we have a number of partial results which may be of interest.  We begin with two rather standard results (perhaps already known, though we could not find them explicitly in the literature).  First, we show that if the solution ever has sufficiently small mass near the origin at a sufficiently late time, then one in fact has scattering and no soliton 
will ever form:

\begin{theorem}\label{easy}  There exist constants $\eps, R > 0$ depending only on the energy $E$
such that the following statement is true: if $u$ is any spherically symmetric solution to \eqref{nls} obeying \eqref{energy-bound}, and
$$ \lim \inf_{t \to +\infty} \int_{B(0,R)} |u(t,x)|^2\ dx \leq \eps,$$
then there exists a function $u_+ \in H^1(\R^3)$ such that
$$ \lim_{t \to +\infty} \| u(t) - e^{it\Delta} u_+(t) \|_{H^1(\R^3)} = 0.$$
Here $B(0,R)$ is the ball $B(0,R) := \{ x \in \R^3: |x| < R \}$.
\end{theorem}

We prove this easy result (a variant of the small data scattering theory, combined with \eqref{radial-sobolev}) in Section \ref{easy-sec}.

From Theorem \ref{easy}, we see that in order to have non-trivial asymptotic behavior (such as the
presence of solitons), we must permanently station some mass near the spatial origin as $t \to +\infty$.  It is then natural to ask what
the asymptotic behavior of this mass is.  One can show without too much difficulty, even without any assumption of spherical symmetry, that the solution decouples into a free solution (which thus decays near the origin,
thanks to local smoothing estimates) and a remainder which is ``weakly bound''
in the sense that it is asymptotically orthogonal\footnote{More precisely, this is an asymptotically bound state at $t = +\infty$.  There is also of course a notion of an asymptotically bound state at $t = -\infty$, and it is an interesting question (which we were unable to address) as to whether these two notions are at all related.} to all free solutions:

\begin{theorem}\label{decoupling}  Let $u$ be any solution to \eqref{nls} obeying 
\eqref{energy-bound} (not necessarily spherically symmetric).  Then there exists a unique decomposition\footnote{This is perhaps not the optimal decomposition of the solution $u$, as it only isolates the linear evolution component of the radiation term, and places everything else (including non-linear self-interactions of the radiation term) into the ``bound state'' term.  However, modulo errors which go to zero in the energy norm as $t \to +\infty$, this seems to be the correct splitting into radiating and non-radiating states; this is confirmed in the radial case by Theorem \ref{main} below.}
\begin{equation}\label{decomposition}
 u(t) = u_\wb(t) + e^{it\Delta} u_+
\end{equation}
into the \emph{weakly bound component} $u_\wb$ and the \emph{linear radiation component} $e^{it\Delta} u_+$,
obeying the energy bounds
\begin{equation}\label{energy-split}
 \| u_+ \|_{H^1(\R^3)}, \sup_{t \in [0,+\infty)} \| u_\wb(t) \|_{H^1(\R^3)} \lesssim 1
\end{equation}
and the asymptotic orthogonality conditions
\begin{equation}\label{asymptotic-l2}
 \lim_{t \to +\infty} \langle u_\wb, e^{it\Delta} f \rangle_{L^2} = 0
\end{equation}
for all $f \in H^{-1}(\R^3)$, or equivalently that
\begin{equation}\label{asymptotic-h1}
 \lim_{t \to +\infty} \langle u_\wb, e^{it\Delta} f \rangle_{H^1} = 0
\end{equation}
for all $f \in H^1(\R^3)$, where
$$ \langle f, g \rangle_{H^1} := \langle f, g \rangle_{L^2} + \langle \nabla f, \nabla g \rangle_{L^2}
= \int_{\R^3} f \overline{g} + \nabla f \cdot \nabla \overline{g}\ dx$$
is the inner product for the Hilbert space $H^1(\R^3)$.  Furthermore, $u_\wb$ is a (local-in-time) approximate
solution to \eqref{nls} in the sense that
\begin{equation}\label{approx-sol}
 \lim_{T \to +\infty} \| (i \partial_t + \Delta) u_\wb - F(u_\wb) \|_{L^1_t H^1_x([T,T+\tau] \times \R^3)} = 0
\end{equation}
for all $\tau > 0$.  Finally, we have the mass decoupling identity
\begin{equation}\label{mass-decoupling}
 MASS(u_+) + \lim_{T \to +\infty} MASS(u_\wb(T)) = MASS(u(0))
\end{equation}
and the energy decoupling identity
\begin{equation}\label{energy-decoupling}
 \frac{1}{2} \int_{\R^3} |\nabla u_+|^2 + \lim_{T \to +\infty} ENERGY(u_\wb(T)) = ENERGY(u(0));
\end{equation}
in particular the above two limits exist.  Finally, if $u$ is spherically symmetric then $u_\wb$
and $u_+$ are also spherically symmetric.
\end{theorem}

We prove this theorem (which is another variant of the small data scattering theory)
in Section \ref{decoupling-sec}.  We remark that \eqref{mass-decoupling}, when combined
with Theorem \ref{easy}, shows that in the spherically symmetric case the asymptotically bound 
mass $\lim_{T \to +\infty} MASS(u_\wb(T))$ is either equal to zero, or is bounded from below\footnote{Strictly speaking, this argument also requires Lemma \ref{energy-local-decay} to dispose of the possible effect of the radiation term $e^{it\Delta} u_+$.  The above remark can also be derived in the non-spherically symmetric case by a standard
perturbation analysis of the free solution $e^{it\Delta} u_+$ at large times, similar to the construction of wave operators for this equation; we omit the details.} by 
some $\eps = \eps(E) > 0$.  As for the asymptotically bound energy, see 
Corollary \ref{energy} below.

Theorem \ref{decoupling} is consistent with the idea of global 
solutions to \eqref{nls} decoupling into a bound state
plus radiation, but the control on the weakly bound state $u_\wb$ given by the above Theorem is very 
unsatisfactory; for instance, we do not know whether $u_\wb$ approaches a \emph{global} 
solution to \eqref{nls} in any reasonable sense (the above result only shows that for any fixed length 
of time $\tau > 0$, the solution $u_\wb$
approaches a solution to \eqref{nls} on intervals $[T, T+\tau]$ for $T$ sufficiently large \emph{depending on} $\tau$).  In particular we do not know whether $u_\wb(t)$ converges asymptotically to a soliton, even if we restrict the time $t$ to a subsequence of times.  However, in the spherically
symmetric case we can prove that $u_\wb$ becomes asymptotically smooth and has some decay:

\begin{theorem}\label{main}  Let $u$ be any spherically symmetric global solution to \eqref{nls} obeying the energy bound \eqref{energy-bound}.  Then for any $J > 1$ and $\delta > 0$, we can decompose the weakly bound component
$u_\wb$ of \eqref{decomposition} further as
\begin{equation}\label{local-decomp}
 u_\wb(t) = u_\bound(t) + o_{\dot H^1(\R^3)}(1)
\end{equation}
for all $t \geq 0$, where $o_{\dot H^1(\R^3)}(1)$ is a time-dependent error whose $\dot H^1(\R^3)$ norm goes to zero as $t \to \infty$,
and $u_\bound(t)$ is a spherically symmetric function (depending on $J$, $\delta$) obeying the symbol-type estimates
\begin{equation}\label{symbol}
 |\nabla_x^j u_\bound(t,x)| \lesssim \langle x \rangle^{-\frac{3}{2} - j + \delta}
\end{equation}
for all $x \in \R^3$, $t \geq T$, and $0 \leq j \leq J$, where $\langle x \rangle := (1 + |x|^2)^{1/2}$ 
and the implicit constants can depend on the exponents $\delta, J$.  Also, we have
\begin{equation}\label{bound-h1}
\sup_{t \to +\infty} \| u_\bound(t) \|_{H^1(\R^3)} \lesssim 1.
\end{equation}
\end{theorem}

This theorem is considerably more technical than the previous two, and is the main result of this paper.  It is proven in Sections \ref{1-sec}, \ref{2-sec}, \ref{3-sec}.  As mentioned earlier, the assumption \eqref{energy-bound} can be
dropped if we mollify the non-linearity smoothly at infinity to decay like $O(|u|^p)$ for $p < 1 + \frac{4}{3}$, 
although this requires a number of technical changes to the argument (e.g. the introduction of the fractional chain 
rule) which we omit.

Note that the error term in \eqref{local-decomp} is only controlled in the homogeneous space\footnote{Of course, since $\dot H^1$ controls $L^6$ by Sobolev embedding, we do get convergence in $H^1(\R^3)$ on any fixed compact set in space, though again this would be more useful if we could improve the bounds on $u_\bound$ in \eqref{symbol} from $O(\langle x \rangle^{-3/2+})$ to
the square-integrable $O(\langle x\rangle^{-3/2-})$.} $\dot H^1(\R^3)$ 
instead of the more natural energy space $H^1(\R^3)$; our arguments have the defect of not being able to control
the very low frequency portion of the solution $u$ in a satisfactory manner (in particular, we could not prevent the very low frequencies of $u_\bound$ from absorbing a non-trivial amount of mass in the asymptotic limit).  This is also 
related to why the bounds
in \eqref{symbol} only give decay of $O(\langle x \rangle^{-3/2+})$ instead of the exponential decay which is to be expected 
by comparison with the soliton solutions.  If one could improve the error control to $H^1$, and also improve
the bounds in \eqref{symbol} from $O(\langle x \rangle^{-3/2+})$ to $O(\langle x \rangle^{-3/2-})$ or better,
then the above theorem would say that the evolution of \eqref{nls} in the bounded energy case \eqref{energy-bound}
has a compact attractor in $H^1$ (once one projects out the linear radiation $e^{it\Delta} u_+$), namely the set generated by \eqref{symbol}.  In principle this makes the
evolution of $u_\bound$ almost periodic, and reduces the
asymptotic study of \eqref{nls} to the dynamical systems problem of understanding the flow on or near this compact
attractor.  It also raises the possibility that one could analyze the asymptotic behavior of the bound portion of
the solution by passing to a weak limit at infinity.  We caution however that the above theorem does not establish
that $u_\bound$ is in any way a global solution to \eqref{nls}; the best we could do was choose $u_\bound$ so that it
obeyed bounds of the form \eqref{approx-sol} (see also Lemma \ref{bc-lemma} below).

Thus Theorem \ref{main} is still far from satisfactory in a number of ways, though it does represent
some progress towards the soliton resolution conjecture.  We speculate that the use of such identities as the
virial identity, pseudoconformal conservation law, or the Morawetz inequalities (none of which are used in the proof of Theorem \ref{main}, which relies almost entirely on an analysis of the high frequencies only) may be useful tools to 
analyze this problem further; we give one example of this in Theorem \ref{pohozaev-thm} below, when we use
virial-type identities to prove an asymptotic Pohozaev identity.  

To put Theorem \ref{main} in some context, we may compare it with the elliptic analogue, which is as follows.  Consider
a solution $u(x)$ to the nonlinear eigenfunction equation \eqref{nl-eigen}
for some $\omega > 0$, where $u$ is assumed to be spherically symmetric and finite energy, but \emph{not} assumed
to be non-negative (since that would force $u$ to be the unique ground state \cite{coff}).  To put it another way, we assume that $u$ is a critical point of the Lagrangian associated to \eqref{nl-eigen} but is not necessarily the minimizer.  The analogue of Theorem \ref{main} would be then that $u$ is smooth and rapidly decreasing.  In the elliptic case, both facts follow easily from rewriting
the equation \eqref{nl-eigen} as $u = (\omega - \Delta)^{-1}( |u|^2 u )$ and then iterating this equation to bootstrap the finite
energy and spherical symmetry assumptions to obtain arbitrary smoothness and decay; see e.g. \cite{caz}.

In the dispersive case, where we have \eqref{nls} instead of \eqref{nl-eigen},
we do not have the luxury of a smoothing operator such as $(\omega - \Delta)^{-1}$.  However, we do have dispersion, which
in principle achieves a similar effect (locally in space, at least) given sufficient amounts of time.  The heuristic
justification for the asymptotic regularity of $u$ is that any high frequency component of $u$ near the origin
should quickly radiate to spatial infinity as $t \to \pm \infty$, but because $u(t)$ has bounded energy as $t \to +\infty$ this radiation should eventually decay in time on each fixed compact region of space (cf. the local smoothing effect \cite{sjolin}, \cite{vega} for the Schr\"odinger equation).  
Thus if we restrict to $t > T$ for $T$ sufficiently large,
the total energy of high frequencies that pass near the origin for times $t > T$ should be small.  The main
difficulty here is to ensure that the high frequency radiation does not interact significantly with the ``bound states''
of the time-dependent Hamiltonian $-\Delta - |u|^2$.  This gets easier to accomplish if one already has some preliminary
regularity result on $u$, which suggests that we can prove the regularity part of Theorem \ref{main} by an iteration
argument starting from \eqref{energy-bound}, similar to how one proceeds in the elliptic case.  It also transpires that the further one gets away from the
origin, the larger the range of frequencies one can classify as ``high frequency'' (roughly speaking, any frequency higher
than $\langle x \rangle^{-1+\delta}$ can escape to spatial infinity either as $t \to +\infty$ or $t \to -\infty$, depending on whether it is outgoing or incoming with respect to the origin).  This explains the
decay factors in \eqref{symbol}.  Unfortunately we do not get infinite decay, because we could not find a dispersive
analogue\footnote{The difficulty is that, unlike the elliptic case, the portion of $u$ near the origin can eventually influence even the very low frequencies $u$ far away from the origin by Duhamel's formula \eqref{duhamel}.  To control this influence it seems that we need to control the spacetime integral of the non-linearity $|u|^2 u$ on large regions of spacetime.  Using the equation \eqref{nls} one can get some reasonable control on such spacetime integrals, but not infinite decay.  Interestingly, in the case of a soliton $u = Q e^{-i\omega t}$ one can get infinite decay on the spacetime integral of $|u|^2 u$ thanks to the time oscillation; if such a phenomenon held for general solutions then we would also get infinite decay of $u_\bound$.} of the elliptic fact that the smoothing operator $(\omega - \Delta)^{-1}$ is highly local.  

As indicated above, the methods used to prove Theorem \ref{main} are fairly elementary, relying principally on the Duhamel
formula, the Strichartz estimates and the dispersive inequality.  There is one somewhat exotic ingredient however,
which is a microlocal decomposition of free spherically symmetric solutions into incoming waves (waves moving towards
the origin) and outgoing waves (waves moving away from the origin).  This decomposition is useful because we have
a freedom when applying the Duhamel formula to evolve either forwards in time (controlling the present from the past)
or backwards in time (controlling the present from the future, taking advantage of the bounded energy hypothesis \eqref{energy-bound} on the future).  To analyze outgoing waves it is easier if one evolves
forwards in time (as the wave becomes incerasingly dispersed), but for incoming waves it is easier to go backwards in time (in order to avoid the wave concentrating at origin, where the non-linear effects are strongest).  To exploit this decomposition we will rely heavily on duality (and in particular
of the unitarity of the free evolution operators $e^{it\Delta}$).  This decomposition may be of independent interest
(for instance, a similar decomposition is utilized in \cite{igor} to establish global-in-time local smoothing and Strichartz estimates on compact non-trapping perturbations of Euclidean space).

One of the main drawbacks of Theorem \ref{main} is the poor control of the decay of $u_\bound$, which is only
shown to be $O(\langle x \rangle^{-3/2+})$. In particular, we were not able to place $x u_\bound$ in $L^2$, which 
seems a natural objective as one may then apply the virial identity or pseudoconformal conservation law 
directly to obtain further information on $u_\bound$.
We were however able to exploit the weaker decay in Theorem \ref{main}, combined with a modified virial identity,
to prove the following statement:

\begin{theorem}[Asymptotic Pohozaev identity]\label{pohozaev-thm} Let $u$ is any spherically symmetric solution to \eqref{nls} obeying \eqref{energy-bound}, and let $u_\bound$ be as in Theorem \ref{main}.  Then we have
$$ \lim_{\tau \to \infty} (\lim \sup_{T \to +\infty} |\frac{1}{\tau} \int_T^{T+\tau} \int_{\R^3} 4 |\nabla u_\bound|^2 - 3 |u_\bound|^4\ dx dt|) = 0.$$
\end{theorem}

As the name of the Theorem suggests, this theorem is the dispersive analogue of the Pohozaev identity
$$ \int_{\R^3} 4 |\nabla u|^2 - 3 |u|^4\ dx = 0$$
for solutions to the non-linear eigenfunction equation \eqref{nl-eigen}, which is proven by multiplying \eqref{nl-eigen}
against $x \cdot \nabla u$ and then integrating by parts.  The above theorem is a (rather weak) generalization
to the dispersive case, saying that the Pohozaev identity is asymptotically verified in an average sense for the bound portion of the solution (a similar estimate holds for $u_\wb$, by \eqref{local-decomp}, Sobolev embedding, and
the fact that both $u_\bound$ and $u_\wb$ are bounded in $H^1$).  The standard
dispersive generalization of the Pohozaev identity is Glassey's virial identity\footnote{The left hand side is also equal to $\partial_{tt} \int |x|^2 |u|^2\ dx$, however in our applications we have nowhere near enough decay to make sense of this quantity and so we will rely instead on the left-hand side as stated in \eqref{virial}.}
\begin{equation}\label{virial}
\partial_t \int 2 x \cdot  \Im(\nabla u \overline{u})\ dx = \int 4 |\nabla u|^2 - 3 |u|^4\ dx,
\end{equation}
which, unsurprisingly, is also obtained by multiplying \eqref{nls} by $x \cdot \nabla u$ and then integrating by parts.
Unfortunately the finite energy assumption on $u$ is not enough to ensure that the integral on the left-hand side
of \eqref{virial} is finite (even with the spherical symmetry assumption), however if we replace $u$ by $u_\bound$
and we add a slight damping weight of $\langle x \rangle^{-\delta}$ then \eqref{symbol} will suffice to make the integral
finite.  We carry out this (standard) procedure and prove Theorem \ref{pohozaev-thm} rigorously in Section \ref{energy-sec}, after a review of virial identities in Section \ref{virial-sec}.  

Combining Theorem \ref{pohozaev-thm} with Theorem \ref{easy} yields the following lower bound on the asymptotically bound energy:

\begin{corollary}\label{energy}  Let $u$ is any spherically symmetric solution to \eqref{nls} 
obeying \eqref{energy-bound}, and such that the asymptotically bound mass $\lim_{t \to +\infty} MASS(u_\wb)$
is non-zero.  Then
there exists $\sigma = \sigma(E) > 0$ such that the asymptotically bound energy is also bounded from zero:
$$ \lim_{t \to +\infty} ENERGY(u_\bound) = \lim_{t \to +\infty} ENERGY(u_\wb) \geq \sigma > 0.$$
\end{corollary}

\begin{proof} Let $u$ be as above.
By Theorem \ref{easy}, there exists $\eps_0 > 0$ and $R > 0$ (depending only on the energy $E$) such that
$$ \lim \inf_{t \to +\infty} \int_{B(0,R)} |u(t,x)|^2\ dx \geq \eps_0.$$
By \eqref{decomposition}, \eqref{local-decomp}, and Lemma \ref{energy-local-decay} below, we thus have
$$ \lim \inf_{t \to +\infty} \int_{B(0,R)} |u_\bound(t,x)|^2\ dx \geq \eps_0.$$
In particular, from H\"older we have
$$ \lim \inf_{t \to +\infty} \int |u_\bound(t,x)|^4\ dx \geq c$$
for some $c = c(\eps_0,R) > 0$.  As a consequence, we have
$$ \lim \inf_{t \to +\infty} \int_{\R^3} 4 |\nabla u_\bound(t)|^2 - 3 |u_\bound(t)|^4\ dx \leq 
\lim \inf_{t \to +\infty} 8 ENERGY(u_\bound) - c.$$
By Theorem \ref{pohozaev-thm} we thus obtain 
$$ \lim \inf_{t \to \infty} 8 ENERGY(u_\bound) - c \geq 0.$$
Observe from \eqref{energy-split}, \eqref{bound-h1} that the error $o_{\dot H^1(\R^3)}(1)$ in \eqref{local-decomp} is 
bounded in $H^1(\R^3)$ norm, and in particular by \eqref{gn} will be decaying in $L^4(\R^3)$ norm.  From this it is easy 
to see that
$$ \lim_{t \to +\infty} ENERGY(u_\bound) - ENERGY(u_\wb) = 0,$$
and the claim then follows from the above estimates (and \eqref{energy-decoupling}, to establish existence of
the limits).
\end{proof}

Corollary \ref{energy} is a variant of the results in \cite{glassey}, \cite{ogawa}, which state that initial data (which is either localized or spherically symmetric) with negative energy will lead to blowup in finite time.  It is also consistent with soliton resolution, since solitons have strictly positive energy in the $L^2$ supercritical case (see e.g. \cite{caz}).  

{\bf Remark.} Except for Theorem \ref{decoupling}, our results are restricted to the spherically symmetric case (although it seems Theorem \ref{easy} generalizes fairly easily, see the remarks at the end of Section \ref{easy-sec}), mainly because \eqref{radial-sobolev} gives us excellent control on where the solution is large (at the origin) and small (everywhere else).  
The general case is of course far more difficult, but one could
speculate that similar results obtain for general bounded energy solutions to \eqref{nls}.  At any given time there should be a bounded number of ``points of concentration'', and the weakly bound component $u_\wb$ should be asymptotically smooth and concentrated near these points, and decay away from these points.  Furthermore these points should move at bounded speeds, though it is not clear to the author whether these speeds will eventually stabilize to be asymptotically constant.
If these points of concentration end up receding from each other, then the behavior of the bound state near each
concentration point should behave much like the analysis given here for the spherically symmetric case 
(after a Gallilean transform), in analogy with the analysis in \cite{schlag}.  Thus ideally one should be able
to reduce the problem of soliton resolution to the case where the only concentration point is the origin, though
as with the results here, to proceed further seems to require some analysis of the dynamical system of the
compact attractor associated to this concentration point, which the author does not know how to pursue.

\section{Notation}

We use the notation $X \lesssim Y$ or $X = O(Y)$ to denote the estimate $X \leq CY$, where $C > 0$ is a constant
that can depend on exponents (such as regularity exponents $j$ or $\alpha$), as well as the 
energy $E$ in \eqref{energy-bound}, but not on other parameters such as time parameters $t$ or functions $u$.

We shall abuse notation and write $f(|x|) = f(x)$ for spherically symmetric functions $f$, i.e. we think of such
functions as living on $[0,+\infty)$ as well as living on $\R^3$.

We fix $\eta$ to denote a smooth spherically symmetric function on $\R^3$ adapted to the ball $B(0,1)$ which equals
$1$ on $B(0,1/2)$.  For any $R > 0$ use $\eta_R$ to denote the rescaling $\eta_R(x) := \eta(x/R)$ of $x$, thus
$\eta_R$ is adapted to $B(0,R)$ and equals 1 on $B(0,R/2)$.  We shall frequently use differences such as
$\eta_R - \eta_r$ or $1-\eta_R$ to localize space smoothly to various annular regions.

We define the Fourier transform $\hat f(\xi)$ on $\R^3$ by
$$ \hat f(\xi) := \int_{\R^3} f(x) e^{-2\pi ix \cdot \xi}\ d\xi.$$
We then define the Littlewood-Paley operators $P_j$ for all integers $j$ by
\begin{equation}\label{lp-def}
\widehat{P_j f}(\xi) := (\eta_{2^{j+1}}(\xi) - \eta_{2^j}(\xi)) \hat f(\xi).
\end{equation}
We record \emph{Bernstein's inequality}
$$ \| f \|_{L^q(\R^3)} \lesssim 2^{(\frac{3}{p}-\frac{3}{q})j} \| f \|_{L^p(\R^3)}$$
whenever $j \in \Z$, $1 \leq p \leq q \leq \infty$, and the Fourier transform of $f$ is supported
on the ball $B(0, C2^j)$ (e.g. $f$ could equal $P_j g$ for some $g$).  This can be proven by applying
Littlewood-Paley type operators and using Young's inequality, as well as the standard bounds for the convolution kernel
of the Littlewood-Paley operators.

We also define the ordinary Sobolev spaces $H^s(\R^3)$ via the norm
$$ \| f \|_{H^s(\R^3)} := \| \langle \xi \rangle^s \hat f \|_{L^2(\R^3)}$$
and the homogeneous spaces $\dot H^s(\R^3)$ by the norm
$$ \| f \|_{\dot H^s(\R^3)} := \| |\xi|^s \hat f \|_{L^2(\R^3)}.$$
We shall unify these two spaces by defining the scaled Sobolev space $H^\alpha_R(\R^3)$
for any spatial scale $R > 0$ by
$$ \| f \|_{H^\alpha_R(\R^3)} := (\int_{\R^3} (|\xi|^2 + \frac{1}{R^2})^\alpha |\hat f(\xi)|^2)^{1/2}.$$
Thus $H^\alpha_1(\R^3)$ is the inhomogeneous space $H^\alpha(\R^3)$, while $\dot H^\alpha(\R^3)$ is in some sense
the limit of $H^\alpha_R(\R^3)$ as $R \to \infty$.  There is also the easily verified scaling relationship
$$ \| f_R \|_{H^\alpha_R(\R^3)} = R^{3/2 - \alpha} \| f \|_{H^\alpha(\R^3)}$$
for all $R$, where $f_R(x) := f(x/R)$.  We remark that $H^\alpha_R$ is the dual of $H^{-\alpha}_R$, and that
we have the estimates
$$ \| f \|_{\dot H^\alpha(\R^3)} \lesssim \| f \|_{H^\alpha_R(\R^3)} \lesssim \| f \|_{H^\alpha(\R^3)} \lesssim 
R^\alpha \| f \|_{H^\alpha_R(\R^3)}$$
when $R \gtrsim 1$ and $\alpha \geq 0$, and dually that
$$ \| f \|_{\dot H^\alpha(\R^3)} \gtrsim \| f \|_{H^\alpha_R(\R^3)} \gtrsim \| f \|_{H^\alpha(\R^3)} \gtrsim 
R^\alpha \| f \|_{H^\alpha_R(\R^3)}$$
when $R \gtrsim 1$ and $\alpha \leq 0$.

We also have the following handy stability lemma with respect to cutoffs.

\begin{lemma}\label{symbol-cutoff}  Suppose that $a(x)$ is a function on $\R^3$ obeying the estimates
$$ \sup_{x \in \R^3} |\nabla^j a(x)| \leq A R^{-j}$$
for some $A \geq 0$ and $R > 0$, and all $0 \leq j \leq J$.  Then we have the estimate
$$ \| af \|_{H^\alpha_R(\R^3)} \lesssim A \| f \|_{H^\alpha_R(\R^3)}$$
for all $-J \leq \alpha \leq J$, where the implicit constants may depend on the exponents $\alpha$ and $J$.
\end{lemma}

\begin{proof}  We may rescale $A=1$ and $R=1$.  By duality we may take $\alpha \geq 0$; by interpolation we may assume that
$0 \leq \alpha \leq J$ is an integer.  But then the claim is clear from the Leibnitz rule and H\"older's inequality.
\end{proof}

We shall also need the standard Sobolev spaces
$$ \| f \|_{W^{k,p}(\R^3)} := \sum_{j=0}^k \| \nabla^j f \|_{L^p(\R^3)}$$
and
$$ \| f \|_{\dot W^{k,p}(\R^3)} := \| \nabla^k f \|_{L^p(\R^3)}$$
for $1 < p < \infty$ and non-negative integers $k$.

\section{Preliminary estimates}

Fix $u$ solving \eqref{nls} and obeying \eqref{energy-bound}; to begin with we do not assume
spherical symmetry.  

Let $e^{it\Delta}$ be the propagator for the free Schr\"odinger equation $iu_t + \Delta u = 0$; this is of
course a unitary operator on $L^2$ (and indeed on every Sobolev space $H^s(\R^3)$).  
From the explicit formula\footnote{We ignore the issue as to which branch of the square root of $t$ to take when $t$ is negative by hiding this issue under the undisclosed constant $C$.}
\begin{equation}\label{explicit}
 e^{it\Delta} f(x) = \frac{C}{t^{3/2}} \int_{\R^3} e^{i|x-y|^2/4t} f(y)\ dy
\end{equation}
we obtain the standard dispersive inequality
\begin{equation}\label{dispersive}
 \| e^{it\Delta} f \|_{L^\infty(\R^3)} \lesssim t^{-3/2} \|f \|_{L^1(\R^3)}
\end{equation}
and as a consequence the global Strichartz inequality (see \cite{tao:keel} and the references therein)\footnote{In this radial setting there are in fact some extra smoothing estimates available,
see \cite{vilela}, although for our argument we will not need any such refinements of Strichartz inequalities (nor will we need more complicated versions involving for instance $X^{s,b}$ spaces.  The main difficulty in
this work is not so much the recovery of local regularity, which is fairly easy since we are in a sub-critical setting with
no derivatives in the nonlinearity, but rather in obtaining sufficient long-term control (e.g. decay) of various components of the solution in time, and the relatively basic Strichartz and dispersive estimates already seem adequate for this task, at least in the three-dimensional case.}
\begin{equation}\label{global-Strichartz-l2}
 \| e^{it\Delta} f \|_{L^2_t L^6_x(\R \times \R^3)} +
 \| e^{it\Delta} f \|_{L^\infty_t L^2_x(\R \times \R^3)}
 \lesssim \| f \|_{L^2(\R^3)},
\end{equation}
and hence (since $e^{it\Delta}$ commutes with derivatives)
\begin{equation}\label{global-Strichartz-h1}
 \| e^{it\Delta} f \|_{L^2_t W^{1,6}_x(\R \times \R^3)} 
+ \| e^{it\Delta} f \|_{L^\infty_t H^1_x(\R \times \R^3)} 
\lesssim \| f \|_{H^1(\R^3)}.
\end{equation}
In particular we observe from Sobolev embedding that
\begin{equation}\label{global-Strichartz-linfty}
 \| e^{it\Delta} f \|_{L^2_t L^\infty_x(\R \times \R^3)} +
 \| e^{it\Delta} f \|_{L^4_t L^6_x(\R \times \R^3)} +
 \| e^{it\Delta} f \|_{L^\infty_t L^3_x(\R \times \R^3)}
\lesssim \| f \|_{H^1(\R^3)}.
\end{equation}
As observed in \cite{tao:keel}, there are inhomogeneous versions of these estimates in which one introduces
a forcing term in a dual Strichartz space.  For instance, we have
\begin{equation}\label{strich-inhomog}
\begin{split}
\| u \|_{L^\infty_t H^1_x([t_0,+\infty) \times \R^3)}& + \| u \|_{L^4_t L^6_x([t_0,+\infty) \times \R^3)}\\
&\lesssim \| u(t_0) \|_{H^1(\R^3)} + \| iu_t + \Delta u \|_{L^2_t W^{1,6/5}_x([t_0,+\infty) \times \R^3)}
\end{split}
\end{equation}
for any (Schwartz) function $u$ and any time $t_0$.  Of course, many other inhomogeneous Strichartz estimates are available, but we isolate this particular one as it shall arise specifically in our arguments\footnote{Note that
these Strichartz estimates are not particularly sharp; we concede a number of derivatives and do not always use the endpoint exponents.  This is a reflection of the fact that the non-linearity is not at a critical power such as
the $L^2$-critical power $p = 1 + \frac{4}{3}$ or the $\dot H^1$-critical power $p=5$.}. 

We also need the following variant of above estimates:

\begin{lemma}\label{4-lemma}
We have
\begin{equation}\label{global-Strichartz-l4}
 \| e^{it\Delta} f \|_{L^4_t L^\infty_x(\R \times \R^3)} \lesssim \| f \|_{\dot H^1(\R^3)}.
\end{equation}
\end{lemma}

\begin{proof} This almost follows by interpolating \eqref{global-Strichartz-h1} with itself, 
but the presence of
the $L^\infty_x$ means that endpoint Sobolev embedding cannot be applied directly, and we must proceed with a bit
more care, using interpolation theory instead.  We shall use the Littlewood-Paley operators $P_j$ defined earlier.  Observe from \eqref{global-Strichartz-l2} and Bernstein's inequality that
\begin{equation}\label{bernstein-2}
 \| e^{it\Delta} P_j f \|_{L^2_t L^\infty_x(\R \times \R^3)} \lesssim 2^{j/2} \| P_j f \|_{L^2(\R^3)}
\sim 2^{-j/2} \| P_j  f\|_{\dot H^1(\R^3)},
\end{equation}
for any $j$, and similarly that
\begin{equation}\label{bernstein-infty}
 \| e^{it\Delta} P_j f \|_{L^\infty_t L^\infty_x(\R \times \R^3)} \lesssim 2^{3j/2} \| P_j f \|_{L^2(\R^3)}
\sim 2^{j/2} \| P_j  f\|_{\dot H^1(\R^3)}.
\end{equation}
By Marcinkeiwicz interpolation we thus obtain the Lorentz space estimate
$$
 \| e^{it\Delta} P_j f \|_{L^{4,1}_t L^\infty_x(\R \times \R^3)} \lesssim \| P_j  f\|_{\dot H^1(\R^3)},
$$ 
and thus by the triangle inequality we have the Besov space estimate
$$
 \| e^{it\Delta} f \|_{L^{4,1}_t L^\infty_x(\R \times \R^3)} \lesssim \sum_{j \in \Z} \| P_j  f\|_{\dot H^1(\R^3)}.
$$ 
On the other hand, from \eqref{bernstein-2} and \eqref{bernstein-infty} we have the estimate
$$ \int_E \| e^{it\Delta} P_j f \|_{L^\infty_x(\R^3)}\ dt \lesssim 
\min(2^{-j/2} |E|^{1/2}, 2^{j/2} |E|) \| P_j  f\|_{\dot H^1(\R^3)}$$
for any measurable set $E$.  Summing this in $j$, we obtain
$$ \int_E \| e^{it\Delta} f \|_{L^\infty_x(\R^3)}\ dt \lesssim 
|E|^{3/4} \sup_j \| P_j  f\|_{\dot H^1(\R^3)}$$
and hence (since $E$ was arbitrary) the weak-type estimate
$$ \| e^{it\Delta} f \|_{L^{4,\infty}_t L^\infty_x(\R^3)}\ dt \lesssim 
\sup_j \| P_j  f\|_{\dot H^1(\R^3)}.$$
Interpolating this with the previous Besov space estimate we obtain the result.
\end{proof}

We can extend these estimates from the linear evolution to the non-linear evolution \eqref{nls}, provided
we localize in time.  The key tool here is the well-known Duhamel formula
\begin{equation}\label{duhamel}
 u(t_1) = e^{i(t_1-t_0) \Delta} u(t_0) - i \int_{t_0}^{t_1} e^{i(t_1-t)\Delta} F(u(t))\ dt
\end{equation}
for solutions to \eqref{nls} and all times $t_0, t_1 \in [0,+\infty)$, where we adopt the convention that
$\int_{t_0}^{t_1} = -\int_{t_1}^{t_0}$ if $t_1 < t_0$.  Of course, it is the time integral in \eqref{duhamel}
which poses the most problems, especially as $|t_1-t_0|$ gets large.  Generally speaking, our strategy in this paper
will be to apply Strichartz estimates or local smoothing estimates in the short term (when $|t_1-t|$ is small)
and to use the dispersive inequality \eqref{dispersive} in the long term (when $|t_1-t|$ is large).

We now give the (standard) non-linear local-in-time analogue of the Strichartz estimates.

\begin{lemma}\label{strichartz-lemma}  Let $u$ be a solution to \eqref{nls} obeying \eqref{energy-bound}.  Then we have
\begin{equation}\label{sobolev-solution}
\sup_{t \in [0,+\infty)} \| u(t) \|_{L^p(\R^3)} \lesssim 1
\end{equation}
for all $2 \leq p \leq 6$, and
\begin{equation}\label{strichartz}
\| u\|_{L^2_t W^{1,6}_x([T,T+\tau] \times \R^3)} \lesssim \langle \tau \rangle^{1/2}
\end{equation}
for all $T \geq 0$ and $\tau > 0$.  In particular we have
\begin{equation}\label{strichartz-linfty}
\| u\|_{L^2_t L^\infty_x([T,T+\tau] \times \R^3)} \lesssim \langle \tau \rangle^{1/2}
\end{equation}
for all $T \geq 0$ and $\tau > 0$.  Similarly, we have
\begin{equation}\label{strichartz-l4}
\| u\|_{L^4_t L^\infty_x([T,T+\tau] \times \R^3)} \lesssim \langle \tau \rangle^{1/4}.
\end{equation}
\end{lemma}

\begin{proof}  The bound \eqref{sobolev-solution} follows immediately from \eqref{energy-bound} and Sobolev.
To prove \eqref{strichartz}, \eqref{strichartz-linfty}, \eqref{strichartz-l4}, it suffices to do so for $\tau > 0$
sufficiently small depending on $E$, since the claim for larger time spans $\tau$ then follows by decomposition of the
time interval into smaller pieces.  We now fix $T$, $\tau$, and let $X$ denote the quantity
$$ X := \| u\|_{L^2_t W^{1,6}_x([T,T+\tau] \times \R^3)} + \| u\|_{L^2_t L^\infty_x([T,T+\tau] \times \R^3)}
+ \| u \|_{L^4_t L^\infty_x([T,T+\tau] \times \R^3)}.$$
By \eqref{global-Strichartz-h1}, \eqref{global-Strichartz-linfty}, \eqref{global-Strichartz-l4}, \eqref{duhamel} and
Minkowski's inequality we have
$$ X \lesssim \| u(T) \|_{H^1(\R^3)} + \int_T^{T+\tau} \| F(u(t)) \|_{H^1(\R^3)}\ dt.$$
By \eqref{energy-bound} the first term on the right-hand side is $O(1)$.  For the second term we use Leibnitz and H\"older to bound
$$ \int_T^{T+\tau} \| F(u(t)) \|_{H^1(\R^3)}\ dt \lesssim \tau^{1/2} \| u \|_{L^4_t L^\infty_x([T,T+\tau] \times \R^3)}^2 \| u \|_{L^\infty_t H^1([T,T+\tau] \times \R^3)}$$
and so by \eqref{energy-bound} and definition of $X$ we have
$$ X \lesssim 1 + \tau^{1/2} X^3$$
and the claim follows by standard continuity arguments (e.g. letting $\tau$ increase continuously from 0) if $\tau$
is sufficiently small.
\end{proof}

Finally, we record a Riemann-Lebesgue type lemma for the free Schr\"odinger propagator $e^{it\Delta}$, which
while preserving the $H^s(\R^3)$ norms, sends many other norms to zero:

\begin{lemma}\label{energy-local-decay}  For any $u_0 \in H^1(\R^3)$ we have
$$ \lim_{t \to +\infty} \int_{\R^3} \langle x \rangle^{-\eps} (|e^{it\Delta} u_0(x)|^2 + |\nabla e^{it\Delta} u_0(x)|^2)\ dx = 0$$
for all $\eps > 0$, and similarly
$$ \lim_{t \to +\infty} \int_{\R^3} |e^{it\Delta} u_0(x)|^p\ dx = 0$$
for all $2 < p \leq 6$.
\end{lemma}

\begin{proof}
Observe that the expressions in the limits are certainly bounded by some quantity depending only
on the $H^1(\R^3)$ norm of $u_0$, thanks to Sobolev embedding and the fact that $e^{it\Delta}$ preserves
the $H^1$ norm.  Thus by the usual limiting argument (and the linearity of $e^{it\Delta}$)
it suffices to verify this lemma for test functions
$u_0$, which are dense in $H^1(\R^3)$.  But in that case we see from \eqref{dispersive} that $|e^{it\Delta} u_0(x)|$
decays like $O(t^{-3/2})$ (with the implicit constants here depending on $u_0$).  
Since $|e^{it\Delta} u_0(x)|$ is also bounded in $L^2$, the claim follows.
\end{proof}

\section{Proof of Theorem \ref{easy}}\label{easy-sec}

We now prove Theorem \ref{easy}, which will be an easy consequence of Strichartz estimates and the dispersive 
inequality.  Let $\eps > 0$
be a small constant to be chosen later, and let $R \gg 1$ be a large number depending on $\eps$ to be chosen later.
From \eqref{global-Strichartz-linfty} we have
$$ \| e^{it\Delta} u(0) \|_{L^4_t L^6_x(\R \times \R^3)} \lesssim 1.$$
Thus (by monotone convergence) we may find a time $T_0 > \eps^{-2} > 0$ (depending on $u$) such that
\begin{equation}\label{linear-ok}
 \| e^{it \Delta} u(0) \|_{L^4_t L^6_x([T_0,+\infty) \times \R^3)} \lesssim \eps.
\end{equation}
By hypothesis, we may then find a time $T_1 \geq T_0$ such that
$$ \int_{B(0,R)} |u(T_1,x)|^2\ dx \lesssim \eps.$$
Then of course we have
\begin{equation}\label{yoink} \int \eta_R |u(T_1,x)|^2\ dx \lesssim \eps,
\end{equation}
where $\eta_R$ was defined in the notation section.  We now bootstrap this control of the local mass
at a fixed time $T_1$ to control of the local mass at nearby times.
From \eqref{nls} we observe the mass flux identity
$$ \partial_t |u|^2 = - 2\nabla \cdot \Im(\nabla u \overline{u})$$
and hence by integration by parts\footnote{One can justify this step by starting with smooth solutions $u$ and then
taking limits using the $H^1$ local well posedness theory.  We will justify similar formal manipulations later in this paper without further comment.}
$$ \partial_t \int \eta_R |u|^2\ dx = 2 \int \nabla \eta_R \cdot \Im(\nabla u \overline{u}).$$
Since $\| \nabla \eta_R \|_{L^\infty} = O(1/R)$, we thus obtain from \eqref{energy-bound} the crude estimate
$$ | \partial_t \int \eta_R |u(t)|^2\ dx| \lesssim 1/R$$
for all times $t$ (far better estimates are available, but this is already sufficient for our purposes).  Thus
by \eqref{yoink} we see that
$$ \sup_{t \in [T_1 - \eps^{-1/4}, T_1]} \int \eta_R |u(t,x)|^2\ dx \lesssim \eps$$
if $R$ is large enough depending on $\eps$.  In particular we see that for any $t \in [T_1 - \eps^{-1/4}, T_1]$ we 
have
$$ \| u(t) \|_{L^3(B(0,R/2))} \lesssim \| u(t) \|_{L^6(\R^3)}^{1/2} \| u(t) \|_{L^2(B(0,R/2))}^{1/2}
\lesssim \eps^{1/2}$$
(using \eqref{sobolev-solution}), while
$$ \| u(t) \|_{L^3(\R^3 \backslash B(0,R/2))} \lesssim \| u(t) \|_{L^\infty(\R^3 \backslash B(0,R/2))}^{2/3} \| u(t) \|_{L^2(\R^3)}^{1/3} \lesssim R^{-2/3}$$
using \eqref{energy-bound}, \eqref{radial-sobolev}.  Thus if $R$ is large enough depending on $\eps$ we have
\begin{equation}\label{l3}
\sup_{t \in [T_1 - \eps^{-1/4}, T_1]}  \| u \|_{L^3(\R^3)} \lesssim \eps^{1/2}.
\end{equation}

We now claim the following smallness bound on the linear evolution (starting from $T_1$).

\begin{lemma}  We have the estimate
\begin{equation}\label{late-linear}
 \| e^{i(t-T_1)\Delta} u(T_1) \|_{L^4_t L^6_x([T_1, +\infty) \times \R^3)}  \lesssim \eps^{1/32}.
\end{equation}
\end{lemma}

\begin{proof}
We first apply Duhamel's formula \eqref{duhamel} to obtain
$$ u(T_1) = e^{iT_1 \Delta} u(0) - i\int_0^{T_1} e^{i(T_1-t')\Delta} F(u(t'))\ dt'$$
and hence
\begin{equation}\label{duhamel-split}
\begin{split}
 e^{i(t-T_1)\Delta} u(T_1) = &e^{it \Delta} u(0)\\
& - i \int_{T_1-\eps^{-1/4}}^{T_1} e^{i(t-t')\Delta} F(u(t'))\ dt'\\
&- i \int_0^{T_1-\eps^{-1/4}} e^{i(t-t')\Delta} F(u(t'))\ dt'.
\end{split}
\end{equation}
From \eqref{linear-ok} we see that the contribution of the first term (the linear component) of \eqref{duhamel-split}
is acceptable.  Now consider the second 
term (the ``recent past'' contribution).  By Strichartz \eqref{global-Strichartz-linfty} and Minkowski's inequality
we see that the contribution of 
this term to \eqref{late-linear} is bounded by
$$ \lesssim \int_{T_1 - \eps^{-1/4}}^{T_1} \| F(u(t')) \|_{H^1(\R^3)}\ dt'.$$
But by the Leibnitz rule and H\"older we can bound this by
$$ \lesssim \| u \|_{L^\infty_t L^3_x} \| u \|_{L^2_t L^\infty_x} \| u \|_{L^2_t W^{1,6}_x}$$
where all norms are over $[T_1 - \eps^{-1/4}, T_1] \times \R^3$.  But by Lemma \ref{strichartz-lemma}  we see that
$$ \| u \|_{L^2_t L^\infty_x}, \| \nabla u \|_{L^2_t L^6_x} \lesssim \eps^{-1/8}$$
and so by combining this with \eqref{l3} we see that this contribution is acceptable.

It remains to prove that
$$ \| \int_0^{T_1-\eps^{-1/4}} e^{i(t-t')\Delta} F(u(t'))\ dt' \|_{L^4_t L^6_x([T_1, +\infty) \times \R^3)}  
\lesssim \eps^{1/32}.$$
On one hand, observe from Duhamel's formula that
$$ \int_0^{T_1-\eps^{-1/4}} e^{i(t-t')\Delta} F(u(t'))\ dt' = e^{i(t-T_1+\eps^{-1/4})\Delta} (u(T_1-\eps^{-1/4}) - u(0))$$
and so from Strichartz \eqref{global-Strichartz-l2}, H\"older, and \eqref{energy-bound} we have
$$ \| \int_0^{T_1-\eps^{-1/4}} e^{i(t-t')\Delta} F(u(t'))\ dt' \|_{L^4_t L^3_x([T_1, +\infty) \times \R^3)}  
\lesssim \| u(T_1-\eps^{-1/4}) - u(0) \|_{L^2(\R^3)} \lesssim 1.$$
Thus by H\"older's inequality again, it will suffice to show that
$$ \| \int_0^{T_1-\eps^{-1/4}} e^{i(t-t')\Delta} F(u(t'))\ dt' \|_{L^4_t L^\infty_x([T_1, +\infty) \times \R^3)}  
\lesssim \eps^{1/16}.$$
But from the dispersive inequality \eqref{dispersive} we have
$$ \| \int_0^{T_1-\eps^{-1/4}} e^{i(t-t')\Delta} F(u(t'))\ dt' \|_{L^\infty(\R^3)}
\lesssim \int_0^{T_1-\eps^{-1/4}} (t-t')^{-3/2} \| F(u(t')) \|_{L^1(\R^3)}\ dt'.$$
But by \eqref{sobolev-solution} we see that $\| F(u(t')) \|_{L^1(\R^3)} = O(1)$, and thus
$$ \| \int_0^{T_1-\eps^{-1/4}} e^{i(t-t')\Delta} F(u(t'))\ dt' \|_{L^\infty(\R^3)}
\lesssim (\eps^{-1/4} + (t-T_1))^{-1/2}.$$
The claim follows.
\end{proof}

We can pass from the bounds \eqref{late-linear} on the linear solution to a corresponding bound
for the non-linear solution by standard arguments.  
Indeed, from the Duhamel formula \eqref{duhamel} we have
$$ u(t) = e^{i(t-T_1)\Delta} u(T_1) - i\int_{T_1}^t e^{i(t-t') \Delta} F(u(t'))\ dt',$$
and upon taking $L^4_t L^6_x([T_1,+\infty) \times \R^3)$ norms, we see from \eqref{late-linear} and
Strichartz \eqref{strich-inhomog} that
$$
 \| u \|_{L^4_t L^6_x([T_1, +\infty) \times \R^3)} \lesssim 
\eps^{1/32} + \| F(u) \|_{L^2_t W^{1,6/5}_x([T_1, +\infty) \times \R^3)}.$$
But by Leibnitz and H\"older we have
$$ \| F(u) \|_{L^2_t W^{1,6/5}_x([T_1, +\infty) \times \R^3)}
\lesssim \| u \|_{L^4_t L^6_x([T_1, +\infty) \times \R^3)}^2
\| u \|_{L^\infty_t H^1_x([T_1, +\infty) \times \R^3)}.$$
Thus from \eqref{energy-bound} and a standard continuity argument (e.g. replacing $[T_1,+\infty)$ by $[T_1,T)$ and then
letting $T \to +\infty$) we thus see that
$$
 \| u \|_{L^4_t L^6_x([T_1, +\infty) \times \R^3)} \lesssim \eps^{1/32}$$
if $\eps$ was chosen sufficiently small.  Finally, we set $u_+$ as
$$ u_+ := e^{-iT_1 \Delta} u(T_1) - i \int_{T_1}^{+\infty} e^{-it' \Delta} F(u(t'))\ dt';$$
the Strichartz inequality \eqref{strich-inhomog} shows that the latter integral converges in $H^1$ since
we have just shown $F(u)$ to lie in $L^2_t W^{1,6/5}_x([T_1,+\infty) \times \R^3)$.  Then 
from the Duhamel formula again we see that
$$ u(t) - e^{-it\Delta} u_+ = i \int_t^{+\infty} e^{-it' \Delta} F(u(t'))\ dt'$$
for all $t \geq T_1$, and thus (by the convergence of the integral) we have
$$ \lim_{t \to +\infty} \| u(t) - e^{it\Delta} u_+ \|_{H^1(\R^3)} = 0$$
as desired.
\endprf

{\bf Remark.}  This argument can mostly be extended to the non-radial case, and shows in fact that one has
scattering to a free solution whenever 
$$\lim \inf_{t \to +\infty} \sup_{x_0 \in \R^3} \int_{B(x_0,R)} |u(t,x)|^2\ dx \leq \eps$$
for some specific $\eps, R > 0$ depending on the energy $E$; in other words, one must always have some concentration 
of the $L^2$ mass in time
in order to prevent the solution from scattering to a free solution, although in the non-radial case one does not have much control on the location $x_0 = x_0(t)$ of this concentration (though the sub-critical nature of this equation, combined with finite speed of propagation heuristics, suggest that we can make $x_0(t)$ depend in a Lipschitz manner in time).  We leave the details of this generalization to the reader (similar results also occur in e.g. \cite{nak:scatter}, \cite{borg:book}).

\section{Proof of Theorem \ref{decoupling}}\label{decoupling-sec}

We now prove Theorem \ref{decoupling}, which is another application of Strichartz and dispersive inequalities.
We begin with a standard preliminary estimate, which in some sense asserts the existence of an adjoint wave operator.

\begin{lemma}\label{adjoint}  Let $u$ be any solution to \eqref{nls} obeying \eqref{energy-bound} (not necessarily
spherically symmetric), and let $f \in H^{-1}(\R^3)$.  Then the limit
\begin{equation}\label{limit}
 \lim_{t \to +\infty} \langle u(t), e^{it\Delta} f \rangle_{L^2}
\end{equation}
exists and is bounded in magnitude by $O(\|f\|_{H^{-1}(\R^3)})$.
\end{lemma}

\begin{proof}  
It is clear that
$$ \sup_{t \in [0,+\infty)} |\langle u(t), e^{it\Delta} f \rangle_{L^2}| \lesssim \| f\|_{H^{-1}(\R^3)}$$
just from \eqref{energy-bound} and duality of $H^1(\R^3)$ and $H^{-1}(\R^3)$.  Thus by the linearity of
$\langle u(t), e^{it\Delta} f \rangle_{L^2}$ in $f$ and a standard limiting argument, it will suffice to prove the 
convergence of the limit \eqref{limit} when $f$ is a test function, i.e. a smooth, compactly supported function.
In particular $f$ now lies in $L^1$ and so from \eqref{dispersive} we have the bound
\begin{equation}\label{decay} \| e^{it\Delta} f \|_{L^\infty} \leq C_f t^{-3/2}
\end{equation}
for all $t \geq 0$, where $C_f$ is a constant depending on $f$.

It will suffice to show that \eqref{limit} is a Cauchy sequence.  Accordingly, we pick an $\eps > 0$ and seek to
find a large time $T > 0$ such that
\begin{equation}\label{inner}
 |\langle u(t), e^{it\Delta} f \rangle_{L^2} - \langle u(T), e^{iT\Delta} f \rangle_{L^2}| \lesssim \eps
\end{equation}
for all $t \geq T$.
To show this we use Duhamel's formula \eqref{duhamel} to write
$$ u(t) = e^{i(t-T)\Delta} u(T) - i \int_T^t e^{i(t-t')\Delta} F(u(t'))\ dt'$$
and hence the left-hand side of \eqref{inner} can be rewritten as
$$ |\langle \int_T^t \langle e^{i(t-t')\Delta} F(u(t')), e^{it\Delta} f \rangle_{L^2}\ dt'|.$$
Thus it will suffice to show that
$$ \int_T^{+\infty} |\langle F(u(t')), e^{it'\Delta} f \rangle_{L^2}|\ dt' \lesssim \eps.$$
But from \eqref{sobolev-solution} we see that $\| F(u)\|_{L^1(\R^3)} \lesssim 1$, while from \eqref{dispersive}
we have $\| e^{it'\Delta} f\|_{L^\infty(\R^3)} \lesssim (t')^{-3/2} \|f\|_{L^1(\R^3)}$.  The claim thus follows
from the decay of $(t')^{-3/2}$ if $T$ is taken sufficiently large.
\end{proof}

The expression \eqref{limit} is thus a bounded linear functional (of $f$) on $H^{-1}(\R^3)$, and hence by
duality there exists a unique $u_+ \in H^1(\R^3)$ such that
$$ \lim_{t \to +\infty} \langle u(t), e^{it\Delta} f \rangle_{L^2} = \langle u_+, f \rangle_{L^2}$$
for all $f \in H^{-1}$, and we have the bound
$\| u_+ \|_{H^1(\R^3)} \lesssim 1$.
If we then define $u_\wb$ by
$$ u_\wb(t) := u(t) - e^{it\Delta} u_+$$
then we clearly have \eqref{decomposition} and \eqref{asymptotic-l2} (and hence \eqref{asymptotic-h1}),
thanks to the unitarity of $e^{it\Delta}$.  Also we obtain \eqref{energy-split} from \eqref{energy-bound}.
Note that this argument also shows that $u_+$ (and hence $u_\wb$) are unique.  Also it is clear
from construction that if $u$ is spherically symmetric then $u_+$ and $u_\wb$ are also.

Now we prove \eqref{approx-sol}.  Fix $\tau$.  Since
$$ (i \partial_t + \Delta) u = F(u)$$
and
$$ (i \partial_t + \Delta) e^{it\Delta} u_+ = 0$$
we have
$$  (i \partial_t + \Delta) u_\wb - F(u_\wb) = F(u) - F(u_\wb).$$
Thus it will suffice to show that
$$ \lim_{T \to \infty} \| F(u) - F(u_\wb) \|_{L^1_t H^1_x([T,T+\tau] \times \R^3)} = 0.$$
Let us just show this for the homogeneous norm $\dot H^1_x$, the contribution of the $L^2_x$ norm being similar.
We observe the pointwise estimate
$$ |\nabla(F(u) - F(u_\wb))| \lesssim (|\nabla u_\wb| + |\nabla u|) 
(|u| + |u_\wb|) |u - u_\wb| + (|u| + |u_\wb|)^2 |\nabla (u - u_\wb)|$$
and hence by H\"older
\begin{align*} \| F(u) - F(u_\wb) \|_{L^1_t \dot H^1_x}
\lesssim &\| |\nabla u_\wb| + |\nabla u| \|_{L^\infty_t L^2_x}
\| |u| + |u_\wb| \|_{L^2_t L^\infty_x} \| u - u_\wb \|_{L^2_t L^\infty_x}\\
&+ \| |u| + |u_\wb| \|_{L^2_t L^6_x} \| |u| + |u_\wb|\|_{L^\infty_t L^6_x}
\| \nabla (u - u_\wb) \|_{L^2_t L^6_x}
\end{align*}
where all norms are in $[T, T + \tau] \times \R^3$.
From \eqref{global-Strichartz-h1} we have
$$ \| \nabla (u - u_\wb) \|_{L^2_t L^6_x(\R \times \R^3)} \lesssim 1$$
and hence
$$ \lim_{T \to \infty} \| \nabla (u - u_\wb) \|_{L^2_t L^6_x([T, T + \tau] \times \R^3)} = 0.$$
A similar argument using \eqref{global-Strichartz-linfty} gives
$$ \lim_{T \to \infty} \| u - u_\wb \|_{L^2_t L^\infty_x([T, T + \tau] \times \R^3)} = 0.$$
Meanwhile from \eqref{global-Strichartz-h1}, \eqref{global-Strichartz-linfty} and Lemma \ref{strichartz-lemma}
we have
$$ \| |u| + |u_\wb| \|_{L^2_t L^p_x([T, T + \tau] \times \R^3)} \leq \langle \tau \rangle^{1/2}$$
for $p = 6, \infty$.  Combining all these bounds together with \eqref{energy-split} and Sobolev embedding
we thus obtain
$$ \lim_{T \to \infty} \| F(u) - F(u_\wb) \|_{L^1_t \dot H^1_x([T,T+\tau] \times \R^3)} = 0$$
as desired.  The analogous argument for the $L^2_x$ norm is similar (because every bound that we have on derivatives
of $u$ and $u_\wb$, we also have on $u$ and $u_\wb$ themselves) and is omitted.  This proves \eqref{approx-sol}.

We now prove \eqref{mass-decoupling}.  Observe that
$$ \|u(t) \|_{L^2(\R^3)}^2 = \| u_\wb(t) \|_{L^2(\R^3)}^2 + \| e^{it\Delta} u_+ \|_{L^2(\R^3)}^2
+ 2 \Re \langle u_\wb(t), e^{it\Delta} u_+ \rangle_{L^2(\R^3)}.$$
The left-hand side is $MASS(u(0))$ by conservation of mass.  The second term on the right-hand side is $MASS(u_+)$ since
$e^{it\Delta}$ is unitary.  The third term on the right goes to zero as $t \to +\infty$ by \eqref{asymptotic-l2},
and the claim follows.

Now we prove \eqref{energy-decoupling}, which is in the same spirit but a little trickier.  We begin with the $H^1$
version of the above identity,
$$ \|\nabla u(t) \|_{L^2(\R^3)}^2 = \| \nabla u_\wb(t) \|_{L^2(\R^3)}^2 + \| \nabla e^{it\Delta} u_+ \|_{L^2(\R^3)}^2
+ 2 \Re \langle \nabla u_\wb(t), \nabla e^{it\Delta} u_+ \rangle_{L^2(\R^3)}.$$
Observe that the third term on the right still goes to zero by \eqref{asymptotic-l2}
(moving both derivatives over to $e^{it\Delta} u_+$).  The left-hand side is equal to $2 ENERGY(u(t)) +
\frac{1}{2} \int |u(t,x)|^4\ dx$, while the first term on the right is equal to $2 ENERGY(u_\wb(t)) + \frac{1}{2} \int |u_\wb(t,x)|^4\ dx$.  Finally, the second term on the right-hand side is just $\int |\nabla u_+(x)|^2\ dx$.  Putting
all this together we see that \eqref{energy-decoupling} will follow if we can show
$$ \lim_{t \to +\infty} \int_{\R^3} |u(t,x)|^4 - |u_\wb(t,x)|^4\ dx = 0.$$
From \eqref{energy-bound}, \eqref{energy-split} we see that 
$\| u(t) \|_{L^4(\R^3)}, \| u_\wb \|_{L^4(\R^3)} \lesssim 1$.  Thus the claim follows from
the pointwise bound
\begin{align}
 \bigl||u(t,x)|^4 - |u_\wb(t,x)|^4\bigr| &\lesssim (|u(t,x)| + |u_\wb(t,x)|)^3
(|u - u_\wb|(t,x)) \\
&= (|u(t,x)| + |u_\wb(t,x)|)^3
|e^{it\Delta} u_+(x)|
\end{align}
and Lemma \ref{energy-local-decay} (applied with $p=4$).
\endprf

{\bf Remark.}  One can obtain more explicit Duhamel-style formulae for $u_\wb$ and $u_+$, namely
$$ u_\wb(t) = i \int_t^{+\infty} e^{i(t-t')\Delta} F(u(t'))\ dt'$$
and
$$ u_+ = u(0) - i \int_0^{+\infty} e^{-it\Delta} F(u(t))\ dt$$
although these integrals do not converge in the energy class, but only in weaker senses such as in the distributional
sense (which is what is essentially done above).  Note that the Duhamel-type formula for $u_\wb$
consists entirely of a forcing term (from the future) and no linear term; we thus expect $u_\wb$ to be
smoother and more localized than $u$ itself, although the infinite time integral $\int_t^{+\infty}$ and the lack of
time decay of $F(u)$ does cause some difficulty in making this heuristic precise.  As indicated in Theorem \ref{main},
we will throw away a small energy error from $u_\wb$ (which is basically caused by the coupling of the
bound state with either the initial data $e^{it\Delta} u(0)$ or the data at infinity $e^{it\Delta} u_+$, as well as the
non-linear self-interactions of the radiation term $e^{it\Delta} u_+$) 
before establishing our final smoothness and decay estimates.

{\bf Remark.}  The above result establishes $L^2$ and $\dot H^1$ type estimates on the radiation term
$u_+$.  It is an interesting problem as to whether there are similar estimates for higher regularities; in particular,
if the map $u(0) \mapsto u_+$ is bounded from $H^s(\R^3)$ to $H^s(\R^3)$ for some $s > 1$.  This type of result is true
in the defocusing case (see e.g. \cite{ckstt:scatter}), but our methods do not seem strong enough to establish this result 
for the focusing equation.  The obstruction is a scenario in which the bound state ($u_\wb$ or $u_\bound$) from continually emits minute amounts of high frequency radiation (which could conceivably cause the $H^s$ norm to become unbounded) without contradicting conservation of mass and energy, and without causing the bound state to collapse into the vacuum state.  We do not know how to show that such a scenario does not occur.

\section{Decomposition into incoming and outgoing waves}\label{1-sec}

We now introduce a fundamental concept in our analysis of the free Schr\"odinger propagator $e^{it\Delta}$, namely the decomposition of a spherically symmetric 
function into \emph{incoming} and \emph{outgoing} components (plus a smooth error); these basically correspond to superpositions of radial waves $e^{i\omega |x|}$ with $\omega < 0$ and $\omega > 0$ 
respectively.   A similar decomposition (constructed using pseudo-differential operators) is also available
for non-spherically-symmetric functions, see for instance \cite{igor}.

\begin{proposition}\label{in-out}  Let $\delta > 0$, and let $R \geq 0$ be such that either $R=0$ or $R \gtrsim 1$.
Let $f$ be a spherically
symmetric test function supported on the exterior region $\{ x: |x| \geq R\}$ (so if $R=0$ then $f$ is just
any spherically symmetric test function).  Then there exists a decomposition
\begin{equation}\label{decomp-pm}
f = f_+ + f_- + f_{smooth}
\end{equation}
with the following properties:

\begin{itemize}

\item ($L^2$-boundedness) We have
\begin{equation}\label{h1-bounds}
 \| f_+ \|_{L^2(\R^3)}, \|f_- \|_{L^2(\R^3)}, \| f_{smooth} \|_{L^2(\R^3)} \lesssim \| f \|_{L^2(\R^3)}.
\end{equation}

\item (Smooth error) The function $f_{smooth}$ obeys the estimates
\begin{equation}\label{smooth-l2}
 \| \nabla f_{smooth} \|_{L^2(\R^3)} \lesssim \langle R \rangle^{(\delta-1)(1+\alpha)} \| f \|_{H_{\langle R\rangle}^{-\alpha}(\R^3)}
\end{equation}
for all $\alpha \geq 0$.

\item (Asymptotic decay of incoming waves) For any $u_0 \in L^2(\R^3)$ and any $\eps > 0$, there exists a time $T = T(u_0,\eps,\delta, R) > 0$ such that
\begin{equation}\label{incoming-decay}
 |\langle e^{it\Delta} u_0, f_-\rangle| \lesssim \eps \| f \|_{L^2(\R^3)}
\end{equation}
for all $t \geq T$ and all test functions $f$ on $\{ x: |x| \geq R\}$.  (Informally, this means that as $t \to +\infty$,
the free wave $e^{it\Delta} u_0$ consists asymptotically of purely outgoing radiation).

\item (Exponential decay near origin in the favorable time direction)  We have
\begin{equation}\label{plus-origin}
 \int_{t_0}^{+\infty} \| \nabla e^{i(t-t_0) \Delta} f_+ \|_{L^2(B(0,R/8))} \ dt 
\lesssim \langle R \rangle^{-\beta} \| f \|_{H^{-\alpha}(\R^3)}
\end{equation}
and
\begin{equation}\label{minus-origin}
 \int_{-\infty}^{t_0} \| \nabla e^{i(t-t_0) \Delta} f_- \|_{L^2(B(0,R/8))} \ dt 
\lesssim \langle R \rangle^{-\beta} \| f \|_{H^{-\alpha}(\R^3)}
\end{equation}
for any $\alpha, \beta \geq 0$.  (Thus, outgoing waves radiate away from the origin as $t \to +\infty$, while incoming waves radiate away from the origin as $t \to -\infty$.)

\item (Local smoothing in the favorable time direction)  We have
\begin{equation}\label{plus-dispersion-alpha}
 \int_{t_0}^{+\infty} \| \langle x \rangle^{-2} e^{i(t-t_0) \Delta} f_+ \|_{H^{-\alpha}_{\langle R \rangle}(\R^3)} \ dt 
\lesssim \langle R \rangle^{-1+\delta} \| f \|_{H^{-\alpha-1+\delta}_{\langle R \rangle}(\R^3)}
\end{equation}
and similarly
\begin{equation}\label{minus-dispersion-alpha}
 \int_{-\infty}^{t_0} \| \langle x \rangle^{-2} e^{i(t-t_0) \Delta} f_- \|_{H^{-\alpha}_{\langle R \rangle}(\R^3)} \ dt 
\lesssim \langle R\rangle^{-1+\delta} \| f \|_{H^{-\alpha-1+\delta}_{\langle R \rangle}(\R^3)}
\end{equation}
for any $\alpha \geq 0$.  Here $H^\alpha_R(\R^3)$ is the scaled Sobolev space defined in the notation section.
\end{itemize}

In the above estimates we allow implicit constants to depend on the exponents $\delta$, $\alpha$, $\beta$ 
but not on $R$.
\end{proposition}

\begin{proof}  In what follows the reader may find it helpful to keep in mind the following heuristics: firstly,
that the Littlewood-Paley projection operators $P_j$ only introduce a spatial uncertainty of $O(2^{-j})$, by
the uncertainty principle; and secondly, the dispersion relation indicates that the Schr\"odinger propagator $e^{i(t-t_0)\Delta}$, when acting on functions of frequency $\sim 2^j$ (such as $P_j f$), will move this function by speeds
comparable to $2^j$.  Because all functions here will be spherically symmetric, there are only two directions for propagation: inwards toward the origin, and outwards away from the origin.  To separate these two modes\footnote{It is possible instead to proceed using the full machinery of Fourier integral operators and semi-classical analysis; for instance, $f_\pm$ is basically the projection of $f$ to the phase space region $|\xi| \gtrsim \langle x \rangle^{-1+\delta}$, $\pm \xi \cdot x \geq 0$, while $f_{smooth}$ is basically the projection to the region $|\xi| \lesssim \langle x \rangle^{-1+\delta}$.  However, to keep the paper reasonably self-contained, we have proceeded more explicitly (though perhaps at the cost of brevity), relying mostly on integration by parts, stationary phase heuristics, and the Fourier inversion formula.  The reader who is content to accept this Proposition may in fact skip the technicalities and move on to the proof of Theorem \ref{main} in the next section.} we will use (truncated) radial Riesz projections.

In the argument that follows we assume some familiarity with the heuristics of stationary phase, and in particular
the ability to get arbitrary decay\footnote{We shall refer to this as the \emph{principle of non-stationary phase}.  We also need a Van der Corput lemma, which asserts that $\int_{\R^3} e^{i\phi(x)} \psi(x) = O(\lambda^{-3/2})$ if $\phi$ has exactly one stationary point with $\nabla^2 \phi \geq \lambda$ at that stationary point, $\psi$ is a rescaled bump function,
and $\phi$ obeys the usual smoothness bounds on the domain of $\psi$.  See \cite{stein:large}.} of oscillatory integrals $\int e^{i\phi(x)} \psi(x)\ dx$ via repeated
integration by parts when the oscillation $\nabla \phi$ of the phase exceeds the amount one loses when 
differentiating $\psi$; we refer to \cite{stein:large} for more precise formulations of this heuristic.

Let us first consider the case when $R=0$, and $f$ is supported on $B(0,2)$.  In this case
we set $f_- := f_{smooth} := 0$ and $f_+ := f$.  The bounds \eqref{h1-bounds}, \eqref{incoming-decay}, \eqref{smooth-l2},
\eqref{minus-dispersion-alpha}
are then trivial, while \eqref{plus-origin}, \eqref{minus-origin} are vacuously true (since $B(0,R/8)$ is empty). 
It thus remains only to verify the outgoing local smoothing estimate \eqref{plus-dispersion-alpha}.  We have to show that
$$
 \int_{t_0}^{+\infty} \| \langle x \rangle^{-2} e^{i(t-t_0) \Delta} f \|_{H^{-\alpha}(\R^3)} \ dt 
\lesssim \| f \|_{H^{-\alpha-1+\delta}(\R^3)}.$$
We use Littlewood-Paley operators to split $f = P_{\leq 0} f + \sum_{j > 0} P_j f$, where
$P_{\leq 0} f := \sum_{j \leq 0} P_j f$.  Consider first the contribution of the low frequencies $P_{\leq 0} f$.
Observe from \eqref{dispersive} that
\begin{align*}
\| e^{i(t-t_0)\Delta} P_{\leq 0} f \|_\infty &\lesssim |t-t_0|^{-3/2} \| P_{\leq 0} f \|_{L^1(B(0,1))}\\
&\lesssim |t-t_0|^{-3/2} \| P_{\leq 0} f \|_{L^2}\\
& \lesssim |t-t_0|^{-3/2} \|f\|_{H^{-\alpha-1+\delta}},
\end{align*}
and also by Bernstein's inequality
$$ \| e^{i(t-t_0)\Delta} P_{\leq 0} f \|_\infty \lesssim \| e^{i(t-t_0)\Delta} P_{\leq 0} f \|_{L^2}
= \| P_{\leq 0} f \|_{L^2} \lesssim \| f \|_{H^{-\alpha-1+\delta}}$$
and hence
$$ \| e^{i(t-t_0)\Delta} P_{\leq 0} f \|_\infty \lesssim \langle t - t_0 \rangle^{-3/2}  \| f \|_{H^{-\alpha-1+\delta}}.$$
By H\"older we thus see that
$$ \| \langle x \rangle^{-2} e^{i(t-t_0)\Delta} P_{\leq 0} f \|_{L^2} \lesssim \langle t - t_0 \rangle^{-3/2}  \| f \|_{H^{-\alpha-1+\delta}}$$
and the claim follows by integrating in $t$ (noting that $L^2$ controls $H^{-\alpha}$).  To control the high frequencies, it will suffice to show that
$$ \int_{t_0}^{+\infty} \| \langle x \rangle^{-2} e^{i(t-t_0) \Delta} P_j f \|_{H^{-\alpha}(\R^3)} \ dt 
\lesssim 2^{(-\alpha+1+\delta/2)j} \| P_j f \|_{L^2(\R^3)}
\lesssim 2^{-\delta j/2} \| f \|_{H^{-\alpha+1-\delta}(\R^3)}$$
for all $j \geq 0$, since the claim then follows by summing in $j$.  Fix $j$; we define $\tau = \tau_j := 2^{(-1+\delta/2)j}$ and divide the above
integral into the immediate future $\int_{t_0}^{t_0 + \tau}$ and the later future $\int_{t_0+\tau}^{+\infty}$. For the immediate future, we use Lemma \ref{symbol-cutoff} to discard $\langle x \rangle^{-2}$, and thus bound this contribution by
$$ \lesssim \int_{t_0}^{t_0 + \tau} \| e^{i(t-t_0) \Delta} P_j f \|_{H^{-\alpha}(\R^3)} \ dt,$$
which is then acceptable by definition of $\tau$ since $e^{i(t-t_0)\Delta}$ preserves the $H^{-\alpha}$ norm.  Now consider the later future $\int_{t_0 + \tau}^{+\infty}$.  For this contribution we write $P_j = \tilde P_j P_j$ where
$\tilde P_j$ is a slight enlargment of the Fourier projection operator $P_j$. Now observe that using stationary
phase that the convolution kernel $K(x)$ of $e^{i(t-t_0)\Delta)} \tilde P_j$ is mostly concentrated in the region
$|x| \gtrsim 2^j |t-t_0| \gtrsim 2^{\delta j/2}$ (the latter bound following definition of $\tau$); indeed
from the principle of non-stationary phase we can obtain bounds of the form $|K(x)| \leq (2^j |t-t_0|)^{-M}$ in the region
$|x| \ll 2^j |t-t_0|$.  Thus the contribution of the region of space $|x| \ll 2^j |t-t_0|$ is easily seen to be manageable.
Thus by applying a smooth cutoff we may restrict to the region $|x| \gtrsim 2^j |t-t_0|$.  But then we may invoke
Lemma \ref{symbol-cutoff} again to replace $\langle x\rangle^{-2}$ by $(2^j |t-t_0|)^{-2}$.  Thus we can bound this contribution by
$$ \lesssim \int_{t_0 + \tau}^{+\infty} (2^j |t-t_0|)^{-2} \| e^{i(t-t_0) \Delta} P_j f \|_{H^{-\alpha}(\R^3)} \ dt,$$
which by the boundedness of $e^{i(t-t_0)\Delta}$ is bounded by
$$ 2^{-2j} \tau^{-1} \| P_j f \|_{H^{-\alpha}(\R^3)}$$
which is acceptable by definition of $\tau$.  This concludes the treatment when $f$ is supported on $B(0,2)$.

It thus remains to consider the case $R \gtrsim 1$, since the case $R=0$ then follows from the $R \gtrsim 1$ case and
the localized case just considered.  We can then replace $\langle R\rangle$ by the comparable quantity $R$ throughout.
We then define $f_{smooth}$ as
$$ f_{smooth} := \sum_{2^j \leq R^{\delta-1}} P_j f + \eta_R \sum_{2^j > R^{\delta-1}} P_j f,$$
where $\eta_R$ is defined in the notation section,
and now proceed to verify the bounds \eqref{smooth-l2} on $f_{smooth}$.  
The contribution of $\sum_{2^j \leq R^{\delta-1}} P_j f$ is
clear just from taking Fourier transforms, so it will suffice to verify that
$$
\sum_{2^j < R^{\delta-1}} \| \nabla \eta_R P_j f \|_{L^2(\R^3)} \lesssim R^{(\delta-1)(1+\alpha)} \| f \|_{H_{R}^{-\alpha}(\R^3)}.$$
But $f$ is supported on $\R^3 \backslash B(0,R)$, which is a distance $R/2$ from the support of $\eta_R$, while
the convolution kernel of $P_j$ decays rapidly outside of $B(0,2^{-j})$, as does all of its derivatives, so it
is easy to verify the pointwise estimate
$$ \| \nabla \eta_R P_j f \|_\infty \lesssim (2^j R)^{-M} \| f \|_{H^{-\alpha}(\R^3)}$$
for any $M > 0$.  Summing over all $j$ with $2^j > R^{\delta - 1}$, we obtain the result (if $M$ is sufficiently large
depending on $\delta$, $\alpha$ we may absorb all powers of $R$ that arise on the right-hand side).  This proves
\eqref{smooth-l2}.

It then remains to decompose
$$ (1-\eta_R) \sum_{2^j > R^{\delta-1}} P_j f = f_- + f_+$$
with the desired properties.  We shall do this by decomposing each $(1-\eta_R) P_j f$ separately 
for each $2^j > R^{\delta-1}$,
\begin{equation}\label{decomp-j}
 (1-\eta_R) P_j f = f_{-,j} + f_{+,j}
\end{equation}
and then defining $f_\pm := \sum_{2^j > R^{\delta - 1}} f_{\pm,j}$.

We begin by taking the radial Fourier transform $g_j$ of the odd extension of $(1-\eta_R) P_j f$,
\begin{equation}\label{gj-def}
 g_j(\rho) := \int_{-\infty}^{+\infty} e^{-2\pi i r \rho} \sgn(r) (1-\eta_R(|r|)) P_j f(|r|)\ dr,
\end{equation}
so that $g_j$ is an odd function on $\R$, and observing the Fourier inversion formula
$$ (1-\eta_R) P_j f = \int_{-\infty}^\infty  g_j(\rho)e_\rho\ d\rho$$
on $\R^3$, where $e_\rho$ is the radial function
$$ e_\rho(x) := e^{2\pi i \rho |x|}.$$
We can thus define
\begin{equation}\label{fpmj-def}
 f_{\pm, j} := \int_{\pm [0,+\infty)} g_j(\rho) (1- \eta_{R/2}) e_\rho \ d\rho
\end{equation}
and we thus obtain the decomposition \eqref{decomp-j} (and hence \eqref{decomp-pm}); note that the
factor $(1-\eta_{R/2})$ does not affect \eqref{decomp-j} since it equals 1 on the support of $(1-\eta_R)$.

Now we verify the $L^2$ bounds \eqref{h1-bounds}.  It suffices to show the claim for $f_\pm$.  By \eqref{fpmj-def}, it suffices to show that
\begin{equation}\label{h1-bounds-pm}
 \| \int_{\pm [0,+\infty)} g(\rho) (1- \eta_{R/2}) e_\rho \ d\rho\|_{L^2(\R^3)} \lesssim \| f \|_{L^2(\R^3)}
\end{equation}
where $g := \sum_{2^j \geq R^{\delta-1}} g_j$.  We discard the $(1-\eta_{R/2})$ and use polar co-ordinates
to write the left-hand side as
$$ \lesssim \| r \int_{\pm [0,\infty)} g(\rho) e^{2\pi i r \rho}\ d\rho \|_{L^2_r}$$
Since $g$ is odd, we can integrate by parts to estimate this as
$$ \lesssim \| \int_{\pm [0,\infty)} g'(\rho) e^{2\pi i r \rho}\ d\rho \|_{L^2_r}$$
which by Plancherel is bounded by
$$ \lesssim (\int_\R |g'(\rho)|^2\ d\rho)^{1/2}.$$
Expanding out $g'$ using \eqref{gj-def}, this is
$$ \lesssim (\int_\R |\int_{-\infty}^{+\infty} e^{-2\pi i r \rho} \sgn(r) (1-\eta_R(|r|)) \sum_{2^j > R^{\delta-1}} P_j f(|r|)\ r dr|^2\ d\rho)^{1/2},$$
which by Plancherel again is bounded by
$$ \lesssim (\int_{-\infty}^{+\infty} |(1-\eta_R(|r|)) \sum_{2^j > R^{\delta-1}} P_j f(|r|)|^2\ r^2 dr)^{1/2}.$$
But if we discard the cutoff $(1-\eta_R)$ and undo the polar co-ordinates, this can be estimated as
$$ \| \sum_{2^j > R^{\delta-1}} P_j f \|_{L^2(\R^3)} \lesssim \|f\|_{L^2(\R^3)}$$
as desired.

Now we prove the asymptotic decay estimate \eqref{incoming-decay} for the incoming wave.  
From \eqref{h1-bounds} we already have
$$
 \sup_{t \geq 0} |\langle e^{it\Delta} u_0, f_-\rangle_{L^2(\R^3)}| \lesssim \| f \|_{L^2(\R^3)} \| u_0 \|_{L^2(\R^3)},
$$
so by the usual Riemann-Lebesgue type limiting argument (as in Lemma \ref{energy-local-decay}) it will suffice to prove \eqref{incoming-decay} assuming
that $u_0$ is a (spherically symmetric) test function.  But in that case $e^{it\Delta} u_0$ has the asymptotics
$$ e^{it\Delta} u_0(x) = t^{-3/2} e^{i|x|^2/4t} a(x/t) + o_{L^2(\R^3)}(1)$$
where $a$ is a spherically symmetric Schwartz function (it is essentially the Fourier transform of $u_0$) and
$o_{L^2(\R^3)}(1)$ is an error whose $L^2(\R^3)$ norm goes to 1 as $t \to +\infty$; see e.g. \cite{heron} (or one can work directly from \eqref{explicit}).  Using \eqref{h1-bounds}
again, it thus suffices to show that
$$ t^{-3/2} |\langle e^{i|x|^2/4t} a(x/t), f_- \rangle_{L^2(\R^3)}| \lesssim \eps \|f\|_{L^2(\R^3)}$$
for all $f \in L^2(\R^3)$ supported on $\R^3 \backslash B(0,R)$, and all $t > T = T(a,\eps,R)$.  We expand out $f_-$
using \eqref{fpmj-def} and estimate the left-hand side by
$$ \lesssim t^{-3/2} \int_{-\infty}^0 |g(\rho)| |\langle e^{i|x|^2/4t} a(x/t), (1-\eta_{R/2}) e_\rho \rangle_{L^2(\R^3)}|\ d\rho.$$
We can write this in polar co-ordinates as
$$
 \lesssim t^{-3/2} \int_{-\infty}^0 |g(\rho)| |\int e^{ir^2/4t - i\rho r} r^2 a(r/t) (1-\eta_{R/2}(r))\ dr|\ d\rho,
$$
and then break this up dyadically as
\begin{equation}\label{previous}
 \lesssim \sum_{k=0}^\infty
t^{-3/2} \int_{-\infty}^0 |g(\rho)| |\int e^{ir^2/4t - i\rho r} r^2 a(r/t) (\eta_{2^k R} -\eta_{2^{k-1} R}(r))\ dr|\ d\rho,
\end{equation}
But the phase $ir^2/4t - i \rho r$ oscillates in $r$ with frequency at least $|\rho| + 2^k R/t$ on
the support of $\eta_{2^k R} - \eta_{2^{k-1} R}$ (because $\rho$ is negative), 
and so by the principle of non-stationary phase we have the bounds
$$ |\int e^{ir^2/4t - i\rho r} r^2 a(r/t) (\eta_{2^k R}-\eta_{2^{k-1} R}(r))\ dr| \lesssim 
\langle 2^k |\rho| R + 2^{2k} R^2 / t \rangle^{-M} (2^k R)^3$$
for any $M > 0$ by repeated integration by parts (the function $a(r/t)$ obeys symbol estimates in $r$ regardless of
what $t$ is, because $a$ is Schwartz).  This decays rapidly for $|\rho| \gg R^{-1} 2^{-k}$.  Hence we can estimate \eqref{previous} as
$$ \lesssim \sum_{k=0}^\infty t^{-3/2} \langle 2^{2k} R^2 / t \rangle^{-M} (2^k R)^3 
(R^{-1} 2^{-k})^{1/2} \| g \|_{L^2_\rho}.$$
This decays rapidly for $2^k \gtrsim t^{1/2} R^{-1}$, so we may sum in $k$ to obtain the bound
$$ \lesssim t^{-3/2} (t^{1/2} R^{-1} R)^3 (R^{-1} t^{-1/2} R)^{1/2} \| g \|_{L^2_\rho} = t^{-1/4} \| g \|_{L^2_\rho}.$$
But by \eqref{gj-def} and (one-dimensional) Plancherel we have the somewhat crude estimate
\begin{equation}\label{crude} \| g_j \|_{L^2_\rho} \lesssim \| (1-\eta_R) P_j f \|_{L^2_r}
\lesssim \| P_j f \|_{L^2(\R^3)}
\end{equation}
and more generally
$$ \| g \|_{L^2_\rho} \lesssim \| (1-\eta_R) \sum_{2^j > R^{\delta-1}} P_j f \|_{L^2_r}
\lesssim \| \sum_{2^j > R^{\delta-1}} P_j f \|_{L^2(\R^3)} \lesssim \|f\|_{L^2(\R^3)}$$
and the claim follows for $t$ sufficiently large.

Now we prove \eqref{plus-origin}, which is basically an application of the principle of non-stationary
phase.  It suffices to show that
$$
 \int_{t_0}^{+\infty} \| \nabla e^{i(t-t_0) \Delta} f_{+,j} \|_{L^2(B(0,R/8))} \ dt 
\lesssim (2^j R)^{-M} \| P_j f \|_{L^2(\R^3)}
$$
for every $2^j \geq R^{\delta-1}$ and every $M > 0$, since the claim then follows by taking $M$ large enough depending on $\alpha, \beta, \delta$ and summing in $j$.

Fix $j$.  By time translation invariance we may take $t_0 = 0$.  
We expand the left-hand side using \eqref{fpmj-def} and use Minkowski's inequality to estimate by
$$\lesssim \int_{0}^{+\infty} \int_0^{+\infty} |g_j(\rho)| 
\| \nabla e^{it\Delta} (1- \eta_{R/2}) e_\rho \|_{L^2(B(0,R/8))} \ d\rho dt.$$
It will suffice to prove the estimate when $(1-\eta_{R/2})$ is replaced by $(\eta_R - \eta_{R/2})$, since
the claim will then follow by replacing $R$ by $2^k R$ for $k = 0, 1, \ldots$ and summing the telescoping series.  By \eqref{explicit}, we have
$$ e^{it\Delta} (\eta_R - \eta_{R/2}) e_\rho(x)
= \frac{C}{t^{3/2}} \int_{\R^3} e^{i|x-y|^2/4t} (\eta_R-\eta_{R/2}(y)) e^{2\pi i \rho |y|}\ dy;$$
using polar co-ordinates in $y$ this becomes
$$ \frac{C}{t^{3/2}} e^{i|x|^2/4t} \int_0^\infty \int_0^\pi 
e^{-i|x| r \cos \theta/2t} (\eta_R-\eta_{R/2}(r)) e^{ir^2/4t} e^{2\pi i \rho r}\ \sin \theta d\theta r^2 dr.$$
Substituting $s = \cos \theta$, this becomes
$$ \frac{C}{t^{3/2}} e^{i|x|^2/4t} \int_{-1}^1 \int_0^\infty 
e^{ir^2/4t} e^{-i|x| rs/2t} e^{2\pi i \rho r} (\eta_R-\eta_{R/2}(r))  \ r^2 dr ds.$$

If $x \in B(0,R/8)$, then the phase $e^{ir^2/4t} e^{-i|x| rs/2t} e^{2\pi i \rho r}$ oscillates in $r$ at a rate of at 
least $\gtrsim \rho + R/t$ on the support of $\eta_R - \eta_{R/2}$ (note that this uses the positivity of $\rho$, i.e. the
fact that we are only considering outgoing waves).  
Meanwhile, every derivative applied to $(\eta_R - \eta_{R/2}(r)) r^2$ lowers this quantity by a factor of about $R$.
Thus one can integrate by parts repeatedly
to obtain a bound of $O(t^{-3/2} R^3 \langle \rho R + R^2/t \rangle^{-M'})$ for any $M'$.  Applying $\nabla_x$ introduces
factors of the order of at worst $R/t$ and thus we obtain a bound of the form
$$ |\nabla e^{it\Delta} (\eta_R - \eta_{R/2}) e_\rho(x)| \lesssim t^{-5/2} R^4 \langle \rho R + R^2/t \rangle^{-M}.$$
Thus it remains to show that
$$ \int_0^{+\infty} \int_{0}^{+\infty} t^{-5/2} R^4 \langle \rho R + R^2/t \rangle^{-M'} |g_j(\rho)|\ d\rho dt
\lesssim (2^j R)^{-M} \| P_j f \|_{L^2(\R^3)}.$$
We can evaluate the $t$ integral to estimate the left-hand side as
$$ \int_0^{+\infty} (\frac{R^2}{\langle R \rho \rangle})^{-3/2} R^4 \langle \rho R \rangle^{-M'} |g_j(\rho)|\ d\rho.$$
If $\rho$ was much larger than $1/R$, e.g. $\rho \geq (2^j/R)^{1/2}$, then the claim is now easily proven from
\eqref{crude}, since $2^j \gtrsim R^{\delta-1}$.  So we may restrict our attention to the 
region $\rho \leq (2^j/R)^{1/2}$.  But then we may use \eqref{gj-def} to write
$$ g_j(\rho) = \int_{-\infty}^{+\infty} e^{-2\pi i r \rho} \sgn(r) (1-\eta_R(|r|)) P_j f(|r|)\ dr,$$
which by collapsing the integral to $\int_0^\infty$ and then undoing the polar co-ordinates becomes
$$ g_j(\rho) = C \int_{\R^3} \sin(2\pi |x| \rho) \frac{1 - \eta_R(x)}{|x|^2} P_j f(x)\ dx.$$
We now integrate by parts repeatedly,
integrating $P_j f$ to effectively gain a factor of $2^{-j}$ (by the frequency localization) while differentiating
$\sin(2\pi |x| \rho)$ or $\frac{1 - \eta_R(x)}{|x|^2}$, at worst losing a factor of $(2^j/R)^{1/2}$.  This
will eventually gain us an arbitrary power of $(2^j R)^{-M}$, and then if we estimate everything else crudely (e.g.
using Cauchy-Schwarz and Plancherel to handle the $\rho$ integral) we will eventually show this contribution is
acceptable (recall that $2^j R > R^\delta$ and so we can eventually overcome any polynomial losses in $R$).
This proves \eqref{plus-origin}.  The proof of \eqref{minus-origin} is very similar and is left to the reader
(it is essentially the conjugate of \eqref{plus-origin}),
so we now turn to \eqref{plus-dispersion-alpha} and \eqref{minus-dispersion-alpha}.  
Actually it will suffice to prove \eqref{plus-dispersion-alpha} since \eqref{minus-dispersion-alpha} is
essentially the conjugate of \eqref{plus-dispersion-alpha} (recall that conjugation turns the forward Schr\"odinger evolution into the backward Schr\"odinger evolution).  

We now prove \eqref{plus-dispersion-alpha}.  Again we may use time translation invariance to set $t_0 = 0$.  
First observe that the portion in the ball $B(0,R/8)$ can be easily treated
by the estimate \eqref{plus-origin} just proven, so we may freely insert a cutoff $(1 - \eta_{R/8})$.  By dyadic decomposition in $j$ it will suffice to show that
$$
 \int_{0}^{+\infty} \| (1 - \eta_{R/8}) \langle x \rangle^{-2} 
e^{it \Delta} f_{+,j} \|_{H^{-\alpha}_{R}(\R^3)} \ dt 
\lesssim R^{-1+\delta} \| P_j f \|_{H^{-\alpha-1}_{R}(\R^3)}$$
for all $2^j \geq R^{\delta-1}$.  

Fix $j$. The intuition here is that at time $t$, the wave $f_{+,j}$ will have propagated away from the origin
by a distance of $\sim 2^j t$, and will hence the weight $\langle x \rangle^{-2}$ will have size 
$O(\langle 2^j t + R \rangle^{-2})$.  Integrating this in time gives us $R^{-1} 2^{-j}$, which explains both
the gain of the power of $R$ and the smoothing of one derivative.  This should be contrasted with the standard
Kato local smoothing result \cite{sjolin}, \cite{vega}, which works both forward and backwards in time (for both incoming and outgoing waves) and gives square-integrable bounds in time, but gains only half a derivative and no gain in $R$.

We turn to the details.  We first split the time integral $\int_0^{+\infty}$ into the near future $\int_0^\tau$
and the distant future $\int_\tau^\infty$, where 
\begin{equation}\label{tau-def}
\tau := 2^{-j} R^{1+\delta/2}
\end{equation}
(this can be interpreted as the time it takes for the wave to leave the region $|x| \sim R$, times a safety margin of $R^{\delta/2}$).  We first consider the near future.
For this region we estimate $(1-\eta_{R/8}) \langle x \rangle^{-2}$ by $R^{-2}$ (using Lemma \ref{symbol-cutoff}) and estimate this contribution by
$$ \int_0^\tau R^{-2} \| e^{it\Delta} f_{+,j} \|_{H^{-\alpha}_R(\R^3)}\ dt.$$
Since $e^{it\Delta}$ preserves the $H^{-\alpha}_R$ norm, it will thus suffice by \eqref{tau-def} to show that
\begin{equation}\label{out-there}
 \| f_{+,j} \|_{H^{-\alpha}_{R}(\R^3)} \lesssim \| P_j f \|_{H^{-\alpha}_R(\R^3)}
\end{equation}
(this estimate will also be useful for the treatment of the distant future, below).
We expand this using \eqref{fpmj-def} as
$$ \| \int_{[0,+\infty)} g_j(\rho) (1- \eta_{R/2}) e_\rho \ d\rho\|_{H^{-\alpha}_R(\R^3)} \lesssim 
2^{-\alpha j} \| P_j f \|_{L^2(\R^3)}.
$$
By an easy modification of the proof of \eqref{h1-bounds} we already have that
$$ \| \int_{[0,+\infty)} g_j(\rho) (1- \eta_{R/2}) e_\rho \ d\rho\|_{L^2(\R^3)} \lesssim 
\| P_j f \|_{L^2(\R^3)}.
$$
so by Littlewood-Paley decomposition (and the hypothesis $2^j \gtrsim R^{1-\delta}$) it will suffice to prove the estimates
\begin{equation}\label{lp-first}
 \| \sum_{j' \leq j-10} P_{j'} \int_{[0,+\infty)} g_j(\rho) (1- \eta_{R/2}) e_\rho \ d\rho\|_{L^2(\R^3)} \lesssim 
(2^{j} R)^{-M} \| P_j f \|_{L^2(\R^3)}
\end{equation}
for all $M$; here of course the implicit constant is allowed to depend on $M$.

This estimate is based on two observations; firstly, that the Fourier transform $((1- \eta_{R/2}) e_\rho)$ is concentrated near the sphere $|\xi| \sim \rho$ (and in particular $P_{j'} ((1-\eta_{R/2}) e_\rho)$ is very small if $\rho \gg 2^{j'}$), and secondly that $g_j(\rho)$ decays very quickly if $\rho/2^j$ is small.
We begin by quantifying the first observation.  We write the Fourier transform $(1- \eta_{R/2}) e_\rho)^\wedge(\xi)$ out as a telescoping sum
$$\sum_{k=0}^\infty \int_{\R^3} (\eta_{2^k R}(y) - \eta_{2^{k-1} R}(y)) e^{2\pi i (\rho |y|-\xi \cdot y)} \ dy.$$
If $\rho > |\xi| + \frac{1}{R}$, then the phase oscillates by at least $\rho-|\xi|$, and so by the principle of non-stationary phase we can bound this expression by
$$ \lesssim \sum_{k=0}^\infty O( (2^k R)^3 (\frac{\rho - |\xi|}{2^k R})^{-M} ) = O(R^3 (\frac{\rho - |\xi|}{R})^{-M})$$
for any large $M$.  In particular we see from Plancherel that
$$ \| \sum_{j' \leq j-10} P_{j'}((1-\eta_{R/2}) e_\rho) \|_{L^2(\R^3)} \lesssim ((2^j + \rho) R)^{-M}$$
when $\rho \geq 2^{j-5}$ for any $M$ (note that this decay can absorb any bounded powers of $2^j$ or $R$).
This easily allows us to control the portion of \eqref{lp-first} arising from the region $\rho \geq 2^{j-5}$, using
\eqref{crude}.

Now we consider the contribution where $\rho \leq 2^{j-5}$.  For this contribution we discard the projector
$\sum_{j' \leq j-10} P_{j'}$ and instead establish decay of $g_j(\rho)$.  By \eqref{gj-def}, \eqref{lp-def}
and the Fourier inversion formula we have
$$ g_j(\rho) = \int_{\R^3} (\eta_{2^{j+1}}(\xi) - \eta_{2^j}(\xi)) \hat f(\xi)
\int_{-\infty}^{+\infty} e^{-2\pi i r \rho} \sgn(r) (1-\eta_R(|r|)) e^{2\pi i |r| e_1 \cdot \xi}
\ dr d\xi$$
where $e_1 = (1,0,0)$ (say). We can divide the $r$ integral into positive and negative axes, and then decompose further dyadically (note there is 
no singularity at the origin because of the cutoff $1-\eta_R(|r|)$.  If $\rho \leq 2^{j-5}$, then the phase oscillates at
a rate of at least $\gtrsim 2^j$, because of the support of the Littlewood-Paley multiplier $\eta_{2^{j+1}} - \eta_{2^j}$.  Thus
by integrating by parts as before, we may bound this expression by
$$ |g_j(\rho)| \lesssim \int_{\R^3} |\eta_{2^{j+1}}(\xi) - \eta_{2^j}(\xi)| |\hat f(\xi)| R (2^j R)^{-M}\ d\xi$$
for any $M > 0$.  Using H\"older and Plancherel, and using the $(2^j R)^{-M}$ factor to absorb any powers of $R$ and $2^j$ which appear, we thus obtain the bound
\begin{equation}\label{gj-small}
 |g_j(\rho)| \lesssim (2^j R)^{-M} \| P_j f\|_{L^2(\R^3)}
\end{equation}
for all $\rho \leq 2^{j-5}$.  In fact a similar argument also gives similar bounds for arbitrarily many derivatives
of $g_j(\rho)$ in this region.  By another application of Plancherel (in the radial direction) we thus have
$$ \| \int_0^{2^{j-5}} g_j(\rho) (1-\eta_{R/2}) e_\rho\ d\rho \|_{L^2(\R^3)} \lesssim (2^j R)^{-M} \| P_j f\|_{L^2(\R^3)}$$
(again absorbing any powers of $2^j$ or $R$ which appear into the $(2^j R)^{-M}$ factor, and so this portion of
the $\rho$ integral is also acceptable.  This proves \eqref{lp-first}, and concludes the treatment of the near future.  

Now we consider the distant future $\int_\tau^\infty$; we have to show that
$$ \int_{\tau}^{+\infty} \| (1 - \eta_{R/8}) \langle x \rangle^{-2} 
e^{it \Delta} f_{+,j} \|_{H^{-\alpha}_{R}(\R^3)} \ dt 
\lesssim R^{-1+\delta} \| P_j f \|_{H^{-\alpha-1}_{R}(\R^3)}$$
We divide the time integral further into dyadic blocks; it will suffice to show that
$$ \int_{2^k \tau}^{2^{k+1} \tau} \| (1 - \eta_{R/8}) \langle x \rangle^{-2} 
e^{it \Delta} f_{+,j} \|_{H^{-\alpha}_{R}(\R^3)} \ dt 
\lesssim 2^{-\delta k} R^{-1+\delta} \| P_j f \|_{H^{-\alpha-1}_{R}(\R^3)}$$
for each $k \geq 0$.

Fix $k$.  Recalling the heuristic that the wave $f_{+,j}$ propagates outward at speed $2^j$, we expect
the function $f_{+,j}$ to be concentrated in the region $|\xi| \gtrsim 2^{j+k} \tau \gg 2^k R$.  Accordingly, 
we split the cutoff $(1-\eta_{R/8})$ into the near component $\eta_{2^k R} - \eta_{R/8}$ and the far component
$1 - \eta_{2^k R}$.

Let us deal with the far component first.  We can estimate $(1 - \eta_{2^k R}) \langle x \rangle^{-2}$ by
$(2^k R)^{-2}$ (using Lemma \ref{symbol-cutoff}), and so we can control this contribution by
$$\lesssim \int_{2^k \tau}^{2^{k+1} \tau} (2^k R)^{-2} \| 
e^{it \Delta} f_{+,j} \|_{H^{-\alpha}_{R}(\R^3)} \ dt.$$
Using the boundedness of $e^{it\Delta}$ on $H^{-\alpha}_R$ and \eqref{out-there} (and then \eqref{tau-def}) we can bound this by
$$ \lesssim 2^k \tau (2^k R)^{-2} 2^{-\alpha j} \| P_j f \|_{L^2(\R^3)}
\sim 2^{-k} R^{-1+\delta/2} \| P_j f \|_{H^{-\alpha-1}_R(\R^3)}$$
which is acceptable.  

Now consider the near component, which will turn out to be very small.  We expand out $f_{+,j}$ using \eqref{fpmj-def} to estimate this contribution (using Minkowski's inequality) by
$$ \lesssim \int_{2^k \tau}^{2^{k+1} \tau} \int_{[0,+\infty)} |g_j(\rho)|
\| (\eta_{2^k R} - \eta_{R/8}) \langle x \rangle^{-2} 
 e^{it\Delta} ((1- \eta_{R/2}) e_\rho) \|_{H^{-\alpha}_{R}(\R^3)} \ dt .$$
We crudely estimate the $H^{-\alpha}_R$ norm by the $L^2(\R^3)$ norm, and 
use Cauchy Schwarz,
combined with the observation that $(\eta_{2^k R} - \eta_{R/8}) \langle x \rangle^{-2}$ has bounded $L^2$ norm 
on $B(0,2^k R)$, to estimate the previous by
\begin{equation}\label{inefficient}
 \lesssim \int_{2^k \tau}^{2^{k+1} \tau} \int_{[0,+\infty)} |g_j(\rho)|
\| e^{it\Delta} ((1- \eta_{R/2}) e_\rho) \|_{L^\infty(B(0,2^k R))} \ dt;
\end{equation}
we shall be able to get away with these very inefficient estimates because we will shortly gain a large power of
$(2^j R)^{-M}$ which can absorb these losses.  Let us first deal with those $\rho$ for which $\rho \geq 2^{j-5}$;
as we saw in the treatment of the immediate future, this is the dominant case.  We now compute $e^{it\Delta} ((1- \eta_{R/2}) e_\rho)(x)$
for $x \in L^2(B(0,2^k R))$, using \eqref{explicit}, as
$$ e^{it\Delta} ((1- \eta_{R/2}) e_\rho)(x) = C t^{-3/2} \int_{\R^3} (1-\eta_{R/2}(x)) e^{2\pi i \rho |y| + i|x-y|^2/4t}\ dy;$$
we may decompose this dyadically as
\begin{equation}\label{dyadic-explicit}
 C t^{-3/2} \sum_{m=0}^\infty \int_{\R^3} (\eta_{2^{m} R}-\eta_{2^{m-1} R}(x)) e^{2\pi i \rho |y| + i|x-y|^2/4t}\ dy.
\end{equation}
Observe that the gradient of the phase $2\pi \rho |y| + |x-y|^2/4t$ in $y$ is
\begin{equation}\label{gradient}
 \nabla_y(2\pi \rho |y| + |x-y|^2/4t) = (\rho + \frac{|y|}{2t}) \frac{y}{|y|} - \frac{x}{2t}.
\end{equation}
Since $\rho \geq 2^{j-5}$, $x \in B(0,2^k R)$, $|y| \sim 2^m R$, and $t \sim 2^k \tau \gg 2^k R / 2^j$, we thus see
that \eqref{gradient} has magnitude $\gtrsim \rho + 2^m R$; the point is that $x/2t$ has magnitude at most
$O(R/\tau) \ll O(\rho)$ and thus does not significantly affect the magnitude of the \eqref{gradient}.  Thus by repeated
integration by parts we may estimate 
$$ |e^{it\Delta} ((1- \eta_{R/2}) e_\rho)(x)| \lesssim 
t^{-3/2} \sum_{m=0}^\infty (2^m R)^3 (2^m R \rho)^{-M}$$
for any $M$; in particular we have
$$ |e^{it\Delta} ((1- \eta_{R/2}) e_\rho)(x)|\lesssim  (\rho R)^{-M} 2^{-3k/2}$$
(using the the fact that $\rho \geq 2^{j-5}$ and $2^j \gtrsim R^{\delta-1}$ to absorb all other factors using the $(\rho R)^{-M}$ term).  Combining this with \eqref{crude} we see
that this contribution to \eqref{inefficient} is
$$ \lesssim 2^{-k/2} R^{-M} \| P_j f \|_{H^{-\alpha-1}_R(\R^3)}$$
for any $M \geq 0$, as desired (note that the factor of $(\rho R)^{-M}$ can absorb any losses of $2^j$ or $R$ which appear).

Finally we consider the portion of the low frequencies where $\rho < 2^{j-5}$.  Here we will enjoy unlimited decay in the $2^j$ and $R$ parameters (because of \eqref{gj-small}) so the only issue is to obtain the required decay in the $2^k$ parameter, which basically amounts to ensuring that one obtains the expected $t^{-3/2}$ decay in the $L^\infty$ norm.
Again we expand $e^{it\Delta} ((1- \eta_{R/2}) e_\rho)(x)$ as \eqref{dyadic-explicit}, but divide now into two cases,
depending on whether $m \leq k+10$ or $m > k+10$.  First suppose $m > k+10$; then the gradient \eqref{gradient}
has magnitude at least $\gtrsim |y|/t \sim 2^{m-k} 2^j R^{-\delta/2}$ (using \eqref{tau-def}), because the $x/2t$ term is much smaller than the other terms.  Thus we may integrate by parts repeatedly and estimate this portion of \eqref{dyadic-explicit} in magnitude by
$$ \lesssim  t^{-3/2} \sum_{m > k+10} (2^m R)^3 (2^m R 2^{m-k} 2^j R^{-\delta/2})^{-M}$$
for arbitrary $M$.  Discarding the $2^{m-k}$ factor we thus see this expression is $O(t^{-3/2} 2^{-Mk} (2^j R)^{-M})$
for any $M$; inserting this bound back into \eqref{inefficient} and using \eqref{gj-small} we see that
this contribution is acceptable.

Now suppose that $m \leq k+10$.  Here the gradient \eqref{gradient} might vanish, but it can only do so at one point
(when $y$ is parallel to $x$ and $\rho + |y|/2t = |x|/2t$).  Furthermore we have the double derivative estimate
$$ \nabla_y^2 (2\pi \rho |y| + |x-y|^2/4t) \geq 1/2t$$
where we use $A \geq B$ to denote the statement that $A-B$ is positive definite; this is basically due to the fact that
the function $2\pi \rho |y|$ is (non-strictly) convex, because $\rho$ is positive.  Hence by Van der Corput's lemma
(see e.g. \cite{stein:large})
we obtain a bound of
$$
|\int_{\R^3} (\eta_{2^{m} R}-\eta_{2^{m-1} R}(x)) e^{2\pi i \rho |y|} e^{i|x-y|^2/4t}\ dy| \lesssim
t^{-3/2};$$
summing this over $m \leq k+10$ we thus see that this contribution to \eqref{dyadic-explicit} is at most $O(\langle k \rangle t^{-3/2})$.  (Actually we may easily get rid of the $\langle k \rangle$, but it will not be important in our argument in any event).  If we combine this with \eqref{gj-small} we obtain the desired result.

This proves \eqref{plus-dispersion-alpha}.  As mentioned earlier, the proof of \eqref{minus-dispersion-alpha}
is similar, so the proof of Proposition
\ref{in-out} is now complete.
\end{proof}

\section{Exterior energy decay of $u_\wb$}\label{2-sec}

We now begin the proof of Theorem \ref{main}, which will be done in a number of stages.  In this section
we will establish the following preliminary decay property of $u_\wb$ away from the origin:

\begin{proposition}\label{exterior-decay}  Let the notation and assumptions be as in Theorem \ref{main}, and let
$0 <  \delta < 1$ be arbitrary.  Then for every $R \geq 0$ we have
$$ \lim \sup_{t \to +\infty} \int_{\R^3 \backslash B(0,R)} |\nabla u_\wb(t,x)|^2\ dx \lesssim R^{-2 + 2\delta}.$$
\end{proposition}

Note that this Proposition only gives decay of the \emph{energy} of $u_\wb$ away from the origin; we suspect that the \emph{mass} of $u_\wb$ also decays away from the origin (which would allow us
to upgrade the $o_{\dot H^1}(1)$ error in \eqref{local-decomp} to the more natural $o_{H^1}(1)$), but we do not know
how to do this.  This Proposition is significantly weaker than Theorem \ref{main} but is also somewhat easier to
prove, and the basic strategy used in this argument will be used again when we prove the rest of Theorem \ref{main}.

\begin{proof}
Fix $\delta, R$.  The claim is trivial from \eqref{energy-split} if $R \lesssim 1$, so we may take $R \gg 1$.
Let $T_0 > 0$ be a large
time depending on $R$, $\delta$ to be chosen later.  We need to show that
$$ \int_{\R^3 \backslash B(0,R)} |\nabla u_\wb(t_0,x)|^2\ dx \lesssim R^{-2+2\delta}$$
for all $t_0 \geq T_0$.
By duality, it suffices to show that
\begin{equation}\label{dual}
 |\langle \nabla u_\wb(t_0), f \rangle_{L^2(\R^3)}| \lesssim R^{-1+\delta}
\end{equation}
for all $t_0 \geq T_0$, and all test functions $f$ supported outside of $B(0,R)$ with bounded 
norm in $L^2(\R^3)$.  Since $u_\wb$ is spherically symmetric we may take $f$ to be also.

Fix $f$, $t_0$.  We use Proposition \ref{in-out} 
to decompose $f = f_- + f_+ + f_{smooth}$, and deal with the 
contribution of each term separately.  The term $f_{smooth}$ is easy, in fact we just use integration by
parts and \eqref{energy-split}
to bound
$$ |\langle \nabla u_\wb(t_0), f_{smooth} \rangle_{L^2(\R^3)}| \lesssim \| \nabla f_{smooth} \|_{L^2(\R^3)}$$
which is acceptable by \eqref{smooth-l2}.

Now we consider the contribution of the outgoing wave $f_+$; here the strategy will be to compare $u$ with the
asymptotic solution $e^{it\Delta} u_+$.  By \eqref{decomposition} it suffices to show that
$$ | \langle \nabla u(t_0), f_+ \rangle_{L^2(\R^3)} - \langle \nabla e^{it_0 \Delta} u_+, f_+ \rangle_{L^2(\R^3)}| \lesssim R^{-1+\delta}
.$$
To prove this we first apply \eqref{duhamel} to obtain
$$
 u(t_0) = e^{i(t_0-T) \Delta} u(T) + i \int_{t_0}^T e^{i(t_0-t)\Delta} F(u(t))\ dt
$$
for any $T \geq 0$. Taking inner products with $f_+$, we thus obtain
$$ \langle \nabla u(t_0), f_+ \rangle_{L^2} = \langle \nabla u(T), e^{i(T-t_0) \Delta} f_+ \rangle_{L^2}
+ i \int_{t_0}^T \int \langle F(u(t)), e^{i(t-t_0) \Delta} f_+ \rangle_{L^2}\ dt.$$
As $T \to +\infty$, we have
$$ \lim_{T \to +\infty} \langle \nabla u(T), e^{i(T-t_0) \Delta} f_+ \rangle_{L^2} =
\langle \nabla u_+, e^{-it_0 \Delta} f_+ \rangle_{L^2} = \langle \nabla e^{it_0\Delta} u_+, f_+ \rangle_{L^2}$$
thanks to \eqref{decomposition}, \eqref{asymptotic-l2}.  Thus it will suffice to show that
$$ \int_{t_0}^\infty |\int \langle \nabla F(u(t)), e^{i(t-t_0) \Delta} f_+ \rangle_{L^2}|\ dt \lesssim R^{-1+\delta}.$$
We now observe that
$$ \| \langle x \rangle^2 \nabla F(u(t)) \|_{H^1(\R^3)} \lesssim 1;$$
on the ball $B(0,1)$ this just comes from \eqref{energy-bound}, \eqref{sobolev-solution} and Leibnitz,
while outside of this ball this comes from \eqref{radial-sobolev} and \eqref{energy-bound}.  The claim then
follows from Cauchy-Schwarz and \eqref{plus-dispersion-alpha} (with $\alpha = 0$).

Finally, we consider the incoming wave $f_-$; which is similar except that we compare $u$ to the initial data $e^{it\Delta} u_0$ instead of the asymptotic data.  From \eqref{incoming-decay} we have
$$ |\langle \nabla e^{it_0 \Delta} u_+, f_- \rangle_{L^2(\R^3)}| \lesssim R^{-1+\delta}$$
and
$$ |\langle \nabla e^{it_0 \Delta} u(0), f_- \rangle_{L^2(\R^3)}| \lesssim R^{-1+\delta}$$
for all $t_0 \geq T_0$, if $T_0$ is chosen sufficiently large.  Thus by \eqref{decomposition} it will
suffice to show that
$$ | \langle \nabla u(t_0), f_- \rangle_{L^2(\R^3)} - \langle \nabla e^{it_0 \Delta} u(0), f_- \rangle_{L^2(\R^3)}| \lesssim R^{-1+\delta}.$$
Applying \eqref{duhamel} we thus reduce to showing that
$$ \int_0^{t_0} |\int \langle \nabla F(u(t)), e^{i(t-t_0) \Delta} f_- \rangle_{L^2}|\ dt \lesssim R^{-1+\delta}.$$
But this follows from \eqref{minus-dispersion-alpha} (with $\alpha = 0$) by arguing as in the $f_+$ case.
This concludes the proof of \eqref{dual} and hence of the Proposition.
\end{proof}

\section{Improving the smoothness of $u_\wb$}\label{3-sec}

We now use an iteration argument to bootstrap the asymptotic regularity of $u_\wb$ from $H^1$
to $H^{1+\alpha}$ for arbitrary values of $\alpha$, and also to gain some symbol-type decay estimates on $u_\wb$ away from the origin.  For inductive purposes it is convenient to phrase
this regularity in a dual formulation.  

\begin{theorem}  Let the notation and assumptions be as in Theorem \ref{main}, and let
$\delta > 0$ and $\alpha \geq 0$ be arbitrary.  Then for every $\eps > 0$ there exists
a time $T = T(u,\delta,\alpha,\eps) > 0$ such that one has
\begin{equation}\label{inductive-regularity}
|\langle \nabla u_\wb(t), f \rangle_{L^2}|
\lesssim \langle R \rangle^{-(1+\alpha)(1-\delta)} \| f \|_{H_{\langle R \rangle}^{-\alpha}} + \eps \| f \|_{L^2}
\end{equation}
for all $t \geq T$, all $R \geq 0$ and all test functions $f$ supported on $\R^3 \backslash B(0,R)$.  Here $H_{\langle R\rangle}^{-1-\alpha}$ is the scaled Sobolev space defined in the notation section.
\end{theorem}

At the end of this section we will use this Theorem to construct $u_\bound$ and thus prove Theorem \ref{main}.

\begin{proof}
We fix $\delta > 0$ and induct on $\alpha$ (keeping $\eps, R$ free to vary).    It suffices by interpolation to prove the claim when
$\alpha$ is a non-negative multiple of $\delta$.  When 
$\alpha = 0$ the claim follows immediately from 
\eqref{energy-split} (just using the first term on the right-hand side of \eqref{inductive-regularity}.
Now suppose inductively that $\alpha \geq \delta$, and the claim \eqref{inductive-regularity} has already been
proven for $\alpha$ replaced by $\alpha-\delta$, and all $\eps, R$ (assuming $T$ sufficiently large of course).
The idea is to use local smoothing and dispersive type estimates and 
Duhamel's identity \eqref{duhamel} to boost the regularity of $u_\wb$ by at least $\delta$ of a 
derivative (the point being that $u_\wb$ does not contain
the linear part of the solution $e^{it\Delta} u_+$ as $t \to +\infty$, which of course has no gain in regularity
except in the local sense of Lemma \ref{energy-local-decay}).

Now fix $\eps$, and let $T > 0$ be a large time (depending on all the above parameters) to be chosen later.  
We first obtain some bounds on the non-linearity $F(u)$ away from the origin.

\begin{lemma}\label{far-out} Let $R \geq 1$.  Then there exists a $T_0 > 0$ such that 
\begin{equation}\label{inductive-nonlinearity}
|\langle \nabla ((1-\eta_{R/8}) F(u(t))), f \rangle_{L^2}|
\lesssim R^{-(1+\alpha-\delta)(1-\delta)} \| \langle x \rangle^{-2} f \|_{H_R^{-\alpha+\delta}} 
+ \eps^{100} \| f \|_{L^2}
\end{equation}
for all $\eps > 0$, $t \geq T_0$, and all test functions $f$ on $\R^3$.
\end{lemma}

\begin{proof}  We first verify this in the case where $f$ vanishes on the ball $B(0, \eps^{-200})$.  Then
on the support of $f$ we have $u = O(\eps^{100})$ by \eqref{radial-sobolev}, and hence by \eqref{energy-bound} we have
$$ \| \nabla ((1-\eta_{R/8}) F(u(t))) \|_{L^2(\R^3)} \lesssim \eps^{100}$$
and the claim follows.  Hence we may assume without loss of generality that $f$ is supported in $B(0, \eps^{-200})$.
Indeed by dyadic decomposition it then suffices to show that
\begin{equation}\label{cut}
|\langle \nabla ((\eta_{R/4} - \eta_{R/8}) F(u(t))), f \rangle_{L^2}|
\lesssim R^{-2} R^{-(1+\alpha-\delta)(1-\delta)} \| f \|_{H_R^{-\alpha+\delta}} 
+ \eps^{200} \| f \|_{L^2}
\end{equation}
for all $1 \leq R \lesssim \eps^{-200}$ equal to a power of 2, and all test functions $f$, as the claim then follows by summing using the triangle inequality (and taking $T_0$ to be the supremum of all the individual values of $T_0$ obtained from each of the finite values of $R$).

Fix $R$.  By the inductive hypothesis, we can ensure that
$$ |\langle \nabla u_\wb(t), f \rangle_{L^2}|
\lesssim R^{-(1+\alpha-\delta)(1-\delta)} \| f \|_{H_{R}^{-\alpha+\delta}} + 
\eps^{500} \| f \|_{L^2}$$
for all $t \geq T_0$ and all test functions $f$ on the annulus $B(0,R) \backslash B(0,R/64)$, 
if $T_0$ is chosen large enough.  By duality this implies
that
$$ (\eta_{R} - \eta_{R/32}) \nabla u_\wb = O_{H_R^{\alpha}}(R^{-(1+\alpha-\delta)(1-\delta)})
+ O_{L^2}(\eps^{500})$$
By Lemma \ref{energy-local-decay}
and \eqref{decomposition} we thus have
$$ (\eta_{R} - \eta_{R/32}) \nabla u = O_{H_R^{\alpha-\delta}}(R^{-(1+\alpha-\delta)(1-\delta)})
+ O_{L^2}(\eps^{500}).$$
if $T_0$ is chosen sufficiently large.  By a Poincare inequality argument (adapted to the scale $R$), this implies that
$$ (\eta_{R/2} - \eta_{R/16}) (u - \overline{u}) = 
O_{H_R^{1+\alpha-\delta}}(R^{-(1+\alpha-\delta)(1-\delta)})
+ O_{H^1_R}(\eps^{500}),$$
where 
$$\overline{u} := \frac{\int_{B(0,R/2) \backslash B(0,R/32)} u}{|B(0,R/2) \backslash B(0,R/32)|} $$
is the mean of $u$ on the annulus $B(0,R/2) \backslash B(0,R/32)$.  Since $R = O(\eps^{-200})$, we may
replace the $O_{H^1_R}(\eps^{500})$ error by $O_{H^1}(\eps^{200})$.
Since $u$ has bounded mass we see that $\overline{u} = O(R^{-3/2})$, and it is then easy to check that 
$(\eta_{R/2} - \eta_{R/16}) \overline{u}$ can be absorbed into the $O_{H_R^{1+\alpha-\delta}}(R^{-(1+\alpha-\delta)(1-\delta)})$ term.  Thus
we have
$$ (\eta_{R/2} - \eta_{R/16}) u = O_{H_R^{1+\alpha-\delta}}(R^{-(1+\alpha-\delta)(1-\delta)})
+ O_{H^1}(\eps^{200}).$$
The contribution of the $O_{H^1}(\eps^{200})$ term to $\nabla((\eta_{R/4}-\eta_{R/8}) |u|^2 u)$ can easily
be seen to be $O_{L^2}(\eps^{200})$, since $u$ is bounded in this region thanks to \eqref{radial-sobolev}.  Thus
by Cauchy-Schwarz the net contribution of this error term is $O(\eps^{200} \|f\|_2)$, which is acceptable, and so
we may ignore this error and pretend in fact that
$$ \| (\eta_{R/2} - \eta_{R/16}) u \|_{H_R^{1+\alpha-\delta}} \lesssim R^{-(1+\alpha-\delta)(1-\delta)}.$$
On the other hand, from \eqref{radial-sobolev} we also have
$$ \| (\eta_{R/2} - \eta_{R/16}) u \|_{L^\infty} \lesssim R^{-1}.$$
From an application of the fractional Leibnitz rule (see e.g. \cite{christ}, \cite{kpv}) we thus have
$$ \| (\eta_{R/4} - \eta_{R/8}) |u|^2 u \|_{H_R^{1+\alpha-\delta}} \lesssim R^{-2} R^{-(1+\alpha-\delta)(1-\delta)},$$
and the claim \eqref{cut} follows by duality.
\end{proof}

We can now verify \eqref{inductive-regularity} near the origin.  This region has to be treated separately
because the radial Sobolev bound \eqref{radial-sobolev} is too singular to be useful here, but the basic structure 
of the argument here will also be used later to prove \eqref{inductive-regularity} away from the origin.

\begin{proposition}\label{near-o}  If $T$ is sufficiently large, then we have
$$ |\langle \nabla u_\wb(t), f \rangle_{L^2}| \lesssim \| f \|_{H^{-\alpha}} + \eps \|f\|_{L^2}$$
for all test functions $f$ supported on $B(0,1)$.
\end{proposition}

\begin{proof}
From the splitting \eqref{decomposition} and Duhamel's formula \eqref{duhamel} we have
$$
\langle \nabla u_\wb(t), f \rangle_{L^2}
= \langle \nabla u_\wb(T'), e^{i(T'-t) \Delta} f \rangle_{L^2}
 + i \int_{t}^{T'} \langle \nabla F(u(t')), e^{i(t'-t)\Delta} f\rangle_{L^2} \ dt'.$$
Letting $T' \to +\infty$ and using \eqref{asymptotic-l2}, it thus suffices to show that
\begin{equation}\label{radiate}
 \int_{t}^{+\infty} |\langle \nabla F(u(t')), e^{i(t'-t)\Delta} f\rangle_{L^2}| \ dt'
\lesssim \| f \|_{H^{-\alpha}} + \eps \|f\|_{L^2}.
\end{equation}
Let $\tau := \eps^{-10}$.  We split this integral into the near future $\int_t^{t+\tau}$ and the
distant future $\int_{t+\tau}^\infty$.  To treat the distant future, we observe from \eqref{energy-bound}, Leibnitz,
Sobolev and H\"older that
$$ \| \nabla F(u(t')) \|_{L^1(\R^3)} \lesssim 1$$
and hence by \eqref{dispersive} we have
$$ \langle \nabla F(u(t')), e^{i(t'-t)\Delta} f\rangle_{L^2} \lesssim |t'-t|^{-3/2}
\| f \|_{L^1} \lesssim |t'-t|^{-3/2} \|f\|_{L^2}$$
by the compact support of $f$.  Thus this contribution is bounded by $\tau^{-1/2} \|f\|_{L^2}$, which is
acceptable by the choice of $\tau$.

It then remains to control the near future $\int_t^{t+\tau}$.
We split $F(u)$ into the portion $\eta_1 F(u)$ near the origin, and the portion $(1-\eta_1) F(u)$
away from the origin. Let us first consider the portion away from the origin.  By \eqref{inductive-nonlinearity}
with $R=1$ we can bound this contribution to \eqref{radiate} by
$$ \lesssim \int_{t}^{t + \tau} \| \langle x \rangle^{-2} e^{i(t'-t)\Delta} f \|_{H^{-\alpha+\delta}} 
+ \eps^{100} \| e^{i(t'-t)\Delta} f \|_{L^2}\ dt'.$$
The second term here is just $\tau \eps^{100} \|f\|_{L^2}$ which is acceptable by choice of $\tau$.
The first term can be controlled using \eqref{plus-dispersion-alpha} as $O( \| f \|_{H^{-\alpha}})$ (recall from the proof of Proposition \ref{in-out} that when $f$ is supported on $B(0,1)$ then we can take $f_+ = f$) which is also
acceptable.  This concludes the treatment of the portion away from the origin.

It remains to show that
$$ \int_{t}^{t + \tau} |\langle \nabla (\eta_1 F(u(t'))), e^{i(t'-t)\Delta} f\rangle_{L^2}| \ dt'
\lesssim \| f \|_{H^{-\alpha}} + \eps \|f\|_{L^2}.$$
From the inductive hypothesis \eqref{inductive-regularity}, with $\alpha$ replaced by $\alpha-\delta$, $\eps$
replaced by $\eps^{100}$, and $R$ replaced by $1$, we have by duality that
$$ \nabla u_\wb(t) = O_{H^{\alpha-\delta}}(1) + O_{L^2}(\eps^{100})$$
if $T$ is large enough.  From this and \eqref{energy-split} we obtain in particular the local estimate
$$ \eta_2 u_\wb(t) = O_{H^{1+\alpha-\delta}}(1) + O_{H^1}(\eps^{100})$$
where $\eta_2$ is a spherically symmetric cutoff to $B(0,1)$ which equals 1 on $B(0,1/2)$.
From \eqref{decomposition} and Lemma \ref{energy-local-decay} we thus see that
$$ \eta_2 u_\wb(t) = O_{H^{1+\alpha-\delta}}(1) + O_{H^1}(\eps^{100})$$
if $T$ is chosen large enough.  

Let us use $u_{error}$ to denote the second term $O_{H^1}(\eps^{100})$.  This error term contributes a number
of terms to $\nabla (\eta_1 F(u))$, which roughly look like either $u^2 \nabla u_{error}$ or 
$u \nabla u u_{error}$ (there is also a term $|u|^2 u_{error}$ arising from when the derivative hits the cutoff).  In all cases we see from \eqref{energy-bound} and Sobolev that these expressions are all
$O_{L^{6/5}}(\eps^{100})$.  On the other hand, by \eqref{global-Strichartz-l2} $e^{i(t'-t)\Delta} f = O_{L^2_t L^6_x}(\|f\|_2)$.
Thus by H\"older the total contribution of the error is at most $O(\eps^{100} \tau^{1/2})$, which is acceptable
by choice of $\tau$.  Thus we shall ignore the error term and pretend that
\begin{equation}\label{hypothesis}
 \| \eta_2 u(t) \|_{H^{1+\alpha-\delta}} \lesssim 1.
\end{equation}
Let us now split up the integral $\int_t^{t+\tau}$ again, into the immediate future $\int_t^{t+1}$ and the
medium-term future $\int_{t+1}^{t+\tau}$.  For the medium-term future we use duality to rewrite this
contribution as
$$ \int_{t+1}^{t+\tau} |\langle e^{i(t-t')\Delta} \eta_1 F(u(t')), \nabla f \rangle_{L^2}|\ dt'.$$
Because the kernel of $e^{i(t-t')\Delta}$ is smooth and decays like $(t-t')^{-3/2}$, it is easy to
verify that
$$ \| \eta_2 e^{i(t-t')\Delta} \eta_1 g \|_{C^K} \lesssim (t-t')^{-3/2} \|g\|_1$$
for any $K > 0$.  In particular we may control this expression in $H^\alpha$.  Since $\||u|^2 u(t')\|_1 \lesssim 1$ 
by \eqref{sobolev-solution}, we thus see that this portion of the integral is bounded 
by $\| f \|_{H^{-\alpha}}$ as desired.

Finally, we consider the immediate future.  For this term we use the well-known 
Kato local smoothing estimate \cite{sjolin}, \cite{vega}
$$ \| \eta_2 e^{i(t'-t)\Delta} f \|_{L^2_{t'} H^{-\alpha+1/2}_x([t,t+1] \times \R^3)}
\lesssim \| f\|_{H^{-\alpha}(\R^3)};$$
this allows us to estimate the contribution of the immediate future by
\begin{equation}\label{immediate} \lesssim \| \nabla (\eta_1 F(u(t'))) \|_{L^2_{t'} H^{\alpha-1/2}_x([t,t+1] 
\times \R^3)}.
\end{equation}
We may freely replace $F(u)$ by $F(\eta_2 u)$.
Heuristically, the fractional Leibnitz rule allows us to write
$$ \nabla^{\alpha-1/2} \nabla (\eta_1 |\eta_2 u|^2 \eta_2 u) \approx \eta_1 (\eta_2 u)^2 \nabla^{\alpha+1/2} (\eta_2 u),$$
which heuristically allows us to estimate \eqref{immediate} by
$$ \lesssim \| \eta_2 u \|_{L^4_t L^\infty_x([t,t+1] \times \R^3)}^2 \| \eta_2 u \|_{L^\infty_t H^{\alpha+1/2}_x([t,t+1] \times \R^3)},$$
which is then acceptable by \eqref{strichartz-l4} and \eqref{hypothesis}, if $\delta$ is sufficiently small.
This heuristic can then be justified by a rigorous application of the fractional Leibnitz rule, for instance by decomposing everything into Littlewood-Paley pieces, or by using the Coifman-Meyer theory of paraproducts; we omit
the details.  (Note that the presence of the $L^\infty_x$ norm here is not dangerous, and in any event we have some surplus
regularity to waste in the $H^{\alpha + 1/2}$ component if one is concerned about endpoints).
\end{proof}

In light of the above proposition, it will suffice to verify \eqref{inductive-regularity} for $R \geq 1$.
Accordingly, let us set $R \geq 1$, $t \geq T$, and let $f$ be supported on $\R^3 \backslash B(0,R)$.  Observe from Proposition \ref{exterior-decay} that
$$
|\langle \nabla u_\wb(t), f \rangle_{L^2}|
\lesssim R^{-1+\delta} \| f \|_{L^2}$$
and so we may assume
\begin{equation}\label{r-bound}
1 \leq R \lesssim \eps^{\frac{-1}{1-\delta}}
\end{equation}
since the claim \eqref{inductive-regularity} follows automatically otherwise from the above discussion.  Indeed the
above argument (combined with a smooth partition of unity) also shows that we may now assume that $f$ is 
supported on the annular region  $\{ R \leq |x| \lesssim \eps^{-1/(1-\delta)}\}$, and we will now do so.
We also observe that it suffices to verify \eqref{inductive-regularity} when $R$ is either equal to zero or an
integer power of two, since the intermediate cases then follow automatically.  In particular, by \eqref{r-bound}
there are only a finite number of values of $R$ to consider.  Thus in order to prove \eqref{inductive-regularity}
it suffices to do so for a fixed value of $R$, since one can then let $T$ be the supremum of all such times $T$
obtained for individual values of $R$.

We thus fix $R, t, f$ and now apply Proposition \ref{in-out} to decompose $f = f_- + f_+ + f_{smooth}$.  The 
contribution of $f_{smooth}$
can be estimated by \eqref{energy-split}, \eqref{smooth-l2} as
$$ |\langle \nabla u_\wb(t), f_{smooth} \rangle_{L^2}| \lesssim \| \nabla f_{smooth} \|_{L^2(\R^3)}
\lesssim R^{-(1+\alpha)(1-\delta)} \| f \|_{\dot H^{-\alpha}}$$
which is acceptable.

We now consider the contribution of $f_-$, arguing as in the proof of Proposition \ref{exterior-decay}.  
From \eqref{incoming-decay} we may assume that
$$
 |\langle \nabla e^{it\Delta} u(0), f_- \rangle| \lesssim \eps \| f \|_{L^2(\R^3)}$$
and
$$
 |\langle \nabla e^{it\Delta} u_+, f_- \rangle| \lesssim \eps \| f \|_{L^2(\R^3)}
$$
if $T$ is chosen sufficiently large.  Thus by \eqref{decomposition} it suffices to show that
$$
|\langle \nabla u(t), f_- \rangle_{L^2} - \langle \nabla e^{it\Delta} u(0), f_- \rangle_{L^2}|
\lesssim R^{-(1+\alpha)(1-\delta)} \| f \|_{H_{R}^{-\alpha}} + \eps \| f \|_{L^2}.$$
But by Duhamel's formula \eqref{duhamel} (as in the proof of Proposition \ref{exterior-decay}) we may estimate the left-hand side
here by
$$ \lesssim \int_0^t |\langle \nabla F(u(t')), e^{i(t'-t)\Delta} f_- \rangle_{L^2}|\ dt'.$$
Let $\tau := \eps^{-10}$.  We shall split this time integral $\int_0^t$ into the \emph{distant past} 
$\int_0^{t-\tau}$ and the \emph{recent past} $\int_{t - \tau}^t$.  Note that $t - \tau$ is positive
if $T$ is chosen sufficiently large depending on $\eps$.  The contribution of the distant past is treated similarly
to the corresponding contribution of the distant future in \eqref{radiate}; we pick up an additional factor of
$O( \eps^{-3/2(1-\delta)} )$ because $f$ is now supported on the ball of radius $O(\eps^{-1/(1-\delta)})$ instead
of $O(1)$, but this still does not affect the argument significantly because $\tau$ is so large.

It remains to consider the recent past.  We have to show that
\begin{equation}\label{recently}
 \int_{t - \tau}^t |\langle \nabla F(u(t')), e^{i(t'-t)\Delta} f_- \rangle_{L^2}|\ dt'
\lesssim R^{-(1+\alpha)(1-\delta)} \| f \|_{H_{R}^{-\alpha}} + \eps \| f \|_{L^2}
\end{equation}
Let $\eta_{R/8}$ be a smooth cutoff to $B(0,R/8)$ which equals 1 on $B(0,R/16)$.  
We split $F(u(t'))$ into the part $\eta_{R/8} F(u(t'))$ near the origin and the part $(1-\eta_{R/8}) F(u(t'))$
away from the origin.  Let us first deal with the part near the origin.  For this part we use Cauchy-Schwarz, observing
from \eqref{energy-bound} and Sobolev that $F(u(t'))$ is in $L^2(\R^3)$, to bound this portion by
$$ \int_{t - \tau}^t \| \nabla e^{i(t'-t)\Delta} f_- \|_{L^2(B(0,R/8)}\ dt'.$$
But this component is acceptable by the exponential decay bound \eqref{minus-origin} near the origin (choosing $\alpha$, $\beta$ in \eqref{minus-origin} sufficiently large).  

It remains to deal with contribution of $(1 - \eta_{R/8}) F(u(t'))$ to \eqref{recently}.  By Lemma \ref{far-out}  we see that this contribution
$$ \lesssim  R^{-(1+\alpha-\delta)(1-\delta)} 
\int_{t - \tau}^t  \| (\langle x \rangle^{-2} e^{i(t'-t)\Delta} f_- \|_{H_R^{-\alpha+\delta}} 
+ \eps^{100} \| e^{(t'-t)\Delta} f_- \|_{L^2})\ dt'.$$
The contribution of the $\eps^{100}$ term is at most
$$\lesssim \eps^{100} \tau \| f_- \|_{\dot H^{-1}},$$
which is acceptable by \eqref{h1-bounds} and the definition of $\tau$.  Now consider the main term.  Applying
\eqref{minus-dispersion-alpha} we may bound this by
$$ \lesssim  R^{-(1+\alpha-\delta)(1-\delta)}  R^{-1+\delta}
\| f_- \|_{H_R^{-\alpha}} $$
which is acceptable.  This proves \eqref{recently} and concludes the treatment of $f_-$.  The treatment of $f_+$ is 
very similar, and involves evolving forward in time instead of backward in time (cf. the proof of
Proposition \ref{near-o}), and so uses $u_+$ instead of $u(0)$ but is otherwise almost identical.  (Note also
that the integration parameter $t'$ is now greater than $t$ and hence also greater than $T$, so there is no difficulty
ensuring that $t'$ is large).  This proves \eqref{inductive-regularity}, which closes the induction and completes the proof of the Theorem.
\end{proof}

We now wrap up by using \eqref{inductive-regularity} to prove Theorem \ref{main}.  We first undo the duality
in \eqref{inductive-regularity}.  Let $\alpha > 0$, $\delta > 0$ and $\eps > 0$ be chosen later, and suppose that $t$ is sufficiently large depending on $\alpha, \delta, \eps$; our implicit constants can depend on $\alpha, \delta$.
For the portion of $u_\wb$ near the origin, we observe from
\eqref{inductive-regularity} (replacing $\eps$ by $\eps^{100}$) that
$$ 
|\langle \eta_2 \nabla u_\wb(t), f \rangle_{L^2}|
\lesssim \| f \|_{H^{-\alpha}} + \eps^{100} \| f \|_{L^2}$$
for all test functions $f$ and for any $\alpha > 0$, $\eps > 0$.  By duality, this means that we may decompose
$$ \eta_2 \nabla u_\wb(t) = O_{H^\alpha(\R^3)}(1) + O_{L^2(\R^3)}(\eps^{100}).$$
On the other hand, from \eqref{energy-split} we have that $\eta_2 u_\wb(t) = O_{H^1(\R^3)}(1)$.
A Poincare inequality argument (writing $\eta_{1/2} u_\wb(t)$ as a localized fractional integration operator
of $\eta_2 \nabla u_\wb(t)$ of order -1, plus an infinitely smoothing operator of $\eta_2 u_\wb(t)$) then easily yields that
$$ \eta_{1/2} u_\wb(t) = O_{H^{\alpha+1}(\R^3)}(1) + O_{H^1(\R^3)}(\eps^{100}).$$
In particular, we can find a function $u_{\bound,0}(t)$ for all sufficiently large times $t$ such that
\begin{equation}\label{close-0}
 \| \eta_1 u_\wb(t) - u_{\bound,0} \|_{H^1(\R^3)} \lesssim \eps^{100}
\end{equation}
and
$$\| u_{\bound,0} \|_{H^{\alpha+1}(\R^3)} \lesssim 1.$$
By applying $\eta_2$ if necessary we may ensure that $u_{\bound,0}$ is supported on $B(0,2)$.  From Sobolev
embedding we see that
\begin{equation}\label{symbol-0}
 |\nabla^j u_{\bound,0}(x)| \lesssim 1
\end{equation}
for all $0 \leq j < \alpha+1/2$.

This is enough to construct $u_\bound$ near the origin.  We now start working in dyadic shells outside of the origin.
We first observe from Proposition \ref{exterior-decay} that it will suffice to construct $u_\bound$ on, say, the ball
$B(0, \eps^{-10})$, since the energy of $u_\wb$ outside this ball is asymptotically zero.  Accordingly,
we let $1/4 \leq R \leq \eps^{-10}$ be a power of 2.  From \eqref{inductive-regularity} we have that
$$
|\langle (\eta_{32R} - \eta_{2R}) \nabla u_\wb(t), f \rangle_{L^2}|
\lesssim R^{-(1+\alpha)(1-\delta)} \| f \|_{H_R^{-\alpha}} + \eps^{100} \| f \|_{L^2}$$
for all test functions $f$, if $t$ is large enough (note that there are only a finite number of $R$ involved, so
we can use a single time threshold $t \geq T = T(\eps,u,\alpha,\delta)$ for all of them).  By duality again,
this implies
\begin{equation}\label{32-form}
 (\eta_{32R} - \eta_{2R}) \nabla u_\wb = O_{H_R^\alpha}( R^{-(1+\alpha)(1-\delta)} )
+ O_{L^2}(\eps^{100}).
\end{equation}
We now claim that a Poincare inequality argument gives
\begin{equation}\label{8-form}
 (\eta_{8R} - \eta_{4R}) u_\wb = O_{H_R^{\alpha+1}}( R^{-(1+\alpha)(1-\delta)} )
+ O_{H^1_R}(\eps^{100}).
\end{equation}
This time we shall perform this argument in detail, taking note of the scaling factor of $R$.  We begin by
observing from integration by parts that for any function $f$ on $\R^3$ and any unit vector $\omega$, we have
$$ f(x) = - \int_0^\infty \eta_R(r) \omega \cdot \nabla f(x + r\omega) + \eta'_R(r) f(x + r\omega)\ dr.$$
Averaging this over all $\omega \in S^2$ and then removing the polar co-ordinates, we obtain
$$ f(x) = C \int_{\R^3} \eta(y) \frac{y}{|y|^3} \cdot \nabla f(x + y) + \frac{\eta'_R(y)}{|y|^2} f(x+y)\ dy.$$
Applying this to $f = u_\wb$ when $x$ lies in the support of $\eta_{8R} - \eta_{4R}$, we see in
particular that
\begin{align*}
(\eta_{8R} - \eta_{4R}) u_\wb(x) =&
C (\eta_{8R}(x) - \eta_{4R}(x)) \int_{\R^3} \eta_R(y) \frac{y}{|y|^3} \cdot \nabla 
((\eta_{32R} - \eta_{2R}) u_\wb)(x+y)\ dy \\
&+C (\eta_{8R}(x) - \eta_{4R}(x)) \int_{\R^3} \eta'_R(y) |y|^{-2} ((\eta_{32R} - \eta_{2R}) u_\wb)(x+y)\ dy.
\end{align*}
We now substitute \eqref{32-form} into the first term, and \eqref{energy-split} into the second term.  Observe that convolution with $\eta_R(y) \frac{y}{|y|^3}$
is a smoothing operator of order -1 and maps $H_R^\alpha$ to $H_R^{\alpha+1}$ with constant $O(1)$ (the factor of $R$
can be easily scaled out) and $L^2$ to $H_R^1$, again with a bound of $O(1)$.  Meanwhile, convolution with
$\eta'_R(y) |y|^{-2}$ is a smoothing operator which maps $L^2$ to $H_R^{\alpha+1}$ with constant $O(1)$.  The
claim \eqref{8-form} follows.

From \eqref{8-form} we may find a function $u_{\bound,R}(t)$ for each $1/4 \leq R \leq \eps^{-10}$ and all sufficiently
large times $t$ such that
\begin{equation}\label{close-R}
 \| (\eta_{8R} - \eta_{4R}) u_\wb - u_{\bound,R} \|_{H^1_R(\R^3)} \lesssim \eps^{100}
\end{equation}
and
$$ \| u_{\bound,R} \|_{H_R^{\alpha+1}} \lesssim R^{-(1+\alpha)(1-\delta)}.$$
By applying $\eta_{16R} - \eta_{2R}$ to $u_{\bound,R}$ if necessary we may assume that $u_{\bound,R}$ is supported
in the annulus $B(0,16R) \backslash B(0,R)$.  From Sobolev embedding (rescaled by $R$) we thus see that
\begin{equation}\label{symbol-R}
 |\nabla^j u_{\bound,R}(t,x)| \lesssim R^{-\frac{3}{2} - j} R^{(1+\alpha)\delta}
\end{equation}
for all integers $0 \leq j < \alpha + 1/2$.  If we now define\footnote{Strictly speaking, this only defines
$u_\bound$ for sufficiently large times.  But for bounded times we can set $u_\bound$ arbitrarily to (say) 0.  In practice
we will actually need to take a sequence of values of $\eps$ tending to zero, with corresponding time thresholds $T_\eps > 0$ tending to infinity, and redefine $u_\bound$ according to a new value of $\eps$ every time we cross one of the time thresholds, in order to make the error in \eqref{local-decomp} $o_{\dot H^1}(1)$ rather than merely $O_{\dot H^1}(\eps)$ for any given $\eps$.  We omit the details.}
$$ u_\bound := u_{\bound,0} + \sum_{1/4 \leq R \leq \eps^{-10}} u_{\bound,R}$$
where $R$ ranges over powers of 2, we thus see by adding up \eqref{symbol-0}, \eqref{symbol-R}
that $u_\bound$ obeys the symbol bounds \eqref{symbol},
if $\alpha$ is chosen sufficiently large depending on $J$, and $\delta$ chosen sufficiently small depending on $\alpha$.
Furthermore, by telescoping \eqref{close-0} and \eqref{close-R} we see that
$$ \| \eta_{8R_0} u_\wb(t) - u_\bound(t) \|_{\dot H^1(\R^3)} \lesssim \eps^{100} \log 1/\eps,$$
where $R_0$ is the largest power of two less than or equal to $\eps^{-10}$.  But from Proposition \ref{exterior-decay}
and Hardy's inequality we see that
$$ \| (1-\eta_{8R_0}) u_\wb(t) \|_{\dot H^1(\R^3)} \lesssim \eps$$
if $t$ is sufficiently large.  Thus we obtain 
$$ u_\wb(t) = u_\bound(t) + O_{\dot H^1}(\eps)$$
which will imply \eqref{local-decomp} since $\eps$ can be arbitrarily small (though of course this forces $t$ to be
increasingly large; see previous footnote).  Finally, from \eqref{close-R} we see that
$$ \| (\eta_{8R} - \eta_{4R}) u_\wb - u_{\bound,R} \|_{H^1(\R^3)} \lesssim R \eps^{100}$$
which telescopes (together with \eqref{close-0}) to obtain that
$$ \| \eta_{8R_0} u_\wb - u_\bound \|_{H^1(\R^3)} \lesssim R_0 \eps^{100}.$$
Since $R_0 \lesssim \eps^{-10}$, the claim \eqref{h1-bounds} then follows from \eqref{energy-split}.  Finally,
we can easily ensure that $u_\bound$ is spherically symmetric (either by making sure each step of the argument
preserves spherical symmetry, or else averaging $u_\bound$ over rotations at the end of the argument. This
concludes the proof of Theorem \ref{main}.
\endprf

\section{A digression on virial-type identities}\label{virial-sec}

Before we begin the proof of Theorem \ref{pohozaev-thm}, we pause to review the derivation of
virial identities such as \eqref{virial}.
First observe that if $u$ is any solution to \eqref{nls},
then we have the momentum flux identity
\begin{equation}\label{flux}
\partial_t \Im(u_j \overline{u}) = -2 \partial_k \Re(\overline{u_j} u_k) + \frac{1}{2} \partial_j \Delta(|u|^2)
+ \frac{1}{2} \partial_j |u|^4
\end{equation}
where $j,k$ range over $1,2,3$ with $u_j := \partial_j u = \partial_{x_j} u$ and the usual summation conventions; 
this can be verified by direct calculation (or see \cite{sulem}).  
Integrating this against a suitable gradient $a_j = \partial_j a$ for some explicit weight $a(x)$ to be chosen 
later, we (formally) obtain after some integration by parts
\begin{align*}
\partial_t \int a_j \Im(u_j \overline{u})\ dx &=
\int -2 a_j \partial_k \Re(\overline{u_j} u_k) + \frac{1}{2} a_j \partial_j \Delta(|u|^2) + \frac{1}{2} a_j \partial_j
|u|^4\ dx \\
&= \int 2 a_{jk} \Re(\overline{u_j} u_k) - \frac{1}{2} (\Delta\Delta a) |u|^2 - \frac{1}{2} (\Delta a) |u|^4\ dx.
\end{align*}
Thus for instance if $a(x) := |x|^2$, then we obtain Glassey's virial identity \eqref{virial}.  Unfortunately we do not have quite enough decay on $u_\bound$ to make the left-hand side of \eqref{virial} finite, and so we will instead use a slightly different weight $a$, namely
$a(x) := |x|^2 \langle x \rangle^{-\delta} \eta_R$, where $\delta > 0$ is a small exponent and $R \gg 1$ is a large
radius, and $\eta_R$ is a cutoff to $B(0,R)$ which equals 1 on $B(0,R/2)$.  Observe that
$$
\Delta \Delta a = O(\delta) \langle x \rangle^{-2-\delta} + O(R^{-2-\delta})$$
and
$$ \Delta a = (6 + O(\delta)) \langle x \rangle^{-\delta} + O(R^{-\delta})$$
on the ball $B(0,R)$, while
$$ a_{jk} \Re(\overline{u_j} u_k) = ((1 + O(\delta)) \langle x \rangle^{-\delta} + O(R^{-\delta})) |\nabla u|^2.$$
We thus have
$$
\partial_t \int 2 a_j \Im(u_j \overline{u})\ dx = \int_{B(0,R)} 4 \langle x \rangle^{-\delta} |\nabla u|^2 
- 3 \langle x \rangle^{-\delta} |u|^4\ dx + O( (\delta + R^{-\delta}) \| u \|_{H^1}^2 + \| u \|_{L^4}^4 ).
$$
We of course can estimate $\| u\|_{L^4}$ by $\|u\|_{H^1}$ by Sobolev embedding.
Integrating this in time, and using the fact that $|\nabla a| = O(\langle x \rangle^{1-\delta}$, we obtain
\begin{equation}\label{virial-modified}
\begin{split}
 |\frac{1}{\tau} \int_T^{T+\tau} \int_{B(0,R)} 4 \langle x \rangle^{-\delta} |\nabla u|^2 
&- 3 \langle x \rangle^{-\delta} |u|^4\ dx dt|
\lesssim \\
&\frac{1}{\tau} \sup_{t = T, T+\tau} \int_{B(0,R)} \langle x \rangle^{1-\delta} |u(t)| |\nabla u(t)|\ dx\\
&+ (\delta + R^{-\delta}) \langle \| u\|_{L^\infty_t H^1_x([T, T + \tau] \times \R^3)} \rangle^4.
\end{split}
\end{equation}
This will be our main tool in proving Theorem \ref{pohozaev-thm}, which we do next.

\section{Proof of Theorem \ref{pohozaev-thm}}\label{energy-sec}

We now prove Theorem \ref{pohozaev-thm}.  Let $u$ be as in the theorem, and let $\eps$ be arbitrary.  We need to show
that for $\tau > 0$ sufficiently large depending on $\eps$, and $T > 0$ sufficiently large depending on $\eps, \tau$, that
\begin{equation}\label{pohozaev-gotit}
 |\frac{1}{\tau} \int_T^{T+\tau} \int_{\R^3} 4 |\nabla u_\bound|^2 - 3 |u_\bound|^4\ dx dt| \lesssim \eps.
\end{equation}

Fix $\eps$.  We now choose a number of parameters; it will be important to pay attention to the order in which these parameters are selected.  For reasons which will be apparent later we will need a large radius $R_0 > 0$ depending on $\eps > 0$ to be chosen later.  To apply \eqref{virial-modified} we shall need a small exponent $\delta > 0$ depending on $\eps$, $R_0$ to be chosen later, as well as an even larger radius $R \gg R_0$ depending on $\delta, \eps, R_0$ to be chosen later, and then
we will choose $\tau$ sufficiently large depending on $\eps, \delta, R_0, R$.  Finally we let $T$ be a sufficiently large time (depending on $u$ and all previous parameters) again to be chosen later.  Let 
$u_{\wb,T}$ be the solution to \eqref{nls} with initial data
$u_{\wb,T}(T) = u_\wb(T)$.  Similarly define $u_{\bound,T}$.   We now show that these
solutions will approximate the actual functions $u_\wb$, $u_\bound$ on the time interval $[T, T+\tau]$
if $T$ is large enough.

\begin{lemma}\label{bc-lemma}  If $T$ is sufficiently large (depending on $\tau$), then we have
$$ \| u_\wb(t) - u_{\wb,T}(t) \|_{H^1} \lesssim \eps^{10} R^{-10}$$
and
\begin{equation}\label{bound-close}
 \| u_\bound(t) - u_{\bound,T}(t) \|_{\dot H^1} \lesssim \eps^{10} R^{-10}
\end{equation}
for all $t \in [T,T+\tau]$.  In particular, from the $H^1$ local well-posedness theory we see that $u_{\bound,T}$ and
$u_{\wb,T}$ are well-defined for all $t \in [T,T+\tau]$ (note that the $L^2$ norm of $u_{\bound,T}$ stays
bounded by mass conservation).
\end{lemma}

\begin{proof}  Fix $T, \tau$; all spacetime norms will be on the slab $[T, T + \tau] \times \R^3$.
We will assume \emph{a priori} that $u_\bound$ and $u_\wb$ exist and are in $H^1$ on the entire
time interval $[T, T+\tau]$, this a priori assumption can then be removed by the usual continuity argument, e.g. letting
$\tau$ increase continuously from zero.

We begin with $u_\wb$.  Write $v = u_\wb - u_{\wb,T}$,
then $v(T) = 0$ and we have the equation
$$ iv_t + \Delta v = F(u_\wb) - F(u_\wb-v) + (i\partial_t + \Delta) u_\wb -
F(u_\wb).$$
However, from \eqref{approx-sol} we have that
$$ \| (i\partial_t + \Delta) u_\wb - F(u_\wb)\|_{L^1_t H^1_x} \lesssim \eps'$$
for an $\eps' > 0$ (much smaller than $\eps$) which we will choose later, if $T$ is sufficiently large depending on $\eps'$, $\tau$, $u$.  
We can thus write the above equation as
$$ iv_t + \Delta v = F(u_\wb) - F(u_\wb-v) + O_{L^1_t H^1_x}(\eps').$$
Let $X$ denote the norm
$$ \| v \|_X := \| v \|_{L^\infty_t H^1_x} + \| v \|_{L^4_t L^\infty_x}.$$
From the Strichartz inequality \eqref{global-Strichartz-l4} (and energy estimates and Duhamel's formula) we have
$$ \| v \|_X
\lesssim \| v(T) \|_{H^1_x} + \| F(u_\wb) - F(u_\wb-v) \|_{L^1_t H^1_x} + \eps';$$
we will refrain for the moment from using the fact that $v(T) = 0$.  However, from Leibnitz and H\"older one easily
checks that
$$ \| F(u_\wb) - F(u_\wb-v) \|_{L^2_t H^1_x}
\lesssim \|v \|_X (\|v\|_X + \| u_\wb \|_X)^2$$
(basically, one applies the derivative in the $H^1$ norm to the cubic expression $F(u_\wb) - F(u_\wb-v)$, and whatever the derivative hits is placed in $L^\infty_t L^2_x$, and the other two factors in $L^4_t L^\infty_x$).  Thus by H\"older we have
$$ \| v \|_X
\lesssim \| v(T) \|_{H^1_x} + \tau^{1/2} \| v \|_X (\|v\|_X + \|u_\wb\|_X)^2 + \eps'.$$
Now observe that Lemma \ref{strichartz-lemma} gives bounds on the $X$ norm of $u$, which are also obeyed by the
free solution $e^{it\Delta} u_+$.  Hence by \eqref{decomposition} we 
have $\| u_\wb\|_X \lesssim \langle \tau \rangle^{1/4}$, thus
$$ \| v \|_X
\lesssim \| v(T) \|_{H^1_x} + \tau^{1/2} \| v \|_X (\|v\|_X + \langle \tau \rangle^{1/4})^2 + \eps'.$$
This estimate does not help us directly when $\tau$ is large.  However, if we decompose the interval
$[T, T + \tau]$ into $O(\tau)$ intervals $[t_1, t_1 + \tau']$ of small length $\tau' \ll 1$, and apply this estimate on each interval, we may show (using continuity arguments in the usual manner) a bound of the form
$$ \| v \|_{X([t_1, t_1 + \tau'] \times \R^3)} \lesssim \| v(t_1) \|_{H^1} + \eps'$$
assuming that $\|v(t_1) \|_{H^1}$ was sufficiently small.
In particular we can control the $H^1$ norm of $v(t_1 +\tau')$ in terms of that of $v(t_1)$, plus an error of $\eps'$.
Iterating this we see that (if $\eps'$ is sufficiently small depending on $\tau'$, $\eps$, $R$) we have
$$ \sup_{t \in [T,T+\tau]} \|v(t) \|_{H^1} \lesssim \eps^{10} R^{-10}$$
as desired.

Now we handle $u_\bound$.  From \eqref{local-decomp} and the bounds just proven, it will suffice to show that
$$ \sup_{t \in [T,T+\tau]} \| u_{\bound,T}(t) - u_{\wb,T}(t) \|_{H^1} \lesssim \eps^{10} R^{-10}$$
for $T$ sufficiently large.  Write $w := u_{\bound,T}(t) - u_{\wb,T}(t)$.  By \eqref{local-decomp}
we may ensure that
$$\|w(T)\|_{\dot H^1} \lesssim \eps'$$
for any $\eps' > 0$ to be chosen later, if $T$ is sufficiently large depending on $\eps'$, $u$.  Since $u_{\bound,T}$
and $u_{\wb,T}$ both obey \eqref{nls}, we see that $w$ obeys the equation
$$ iw_t + \Delta w = F(u_{\wb,T}+w) - F(u_{\wb,T}).$$
We now perform a homogeneous version of the previous analysis.  Let $\dot X$ denote the norm
$$ \| w \|_{\dot X} := \| w \|_{L^\infty_t \dot H^1_x} + \| w \|_{L^4_t L^\infty_x},$$
then by Strichartz again (in particular \eqref{global-Strichartz-l4}) we have
$$ \| w \|_{\dot X} \lesssim \| w(T) \|_{\dot H^1} + \| F(u_{\wb,T}+w) - F(u_{\wb,T}) \|_{L^1_t \dot H^1_x}.$$
By Leibnitz and H\"older we have
$$ \| F(u_{\wb,T}+w) - F(u_{\wb,T}) \|_{L^1_t \dot H^1_x} \lesssim
\tau^{1/2} \|w\|_{\dot X} (\| u_{\wb,T} \|_{\dot X} + \| w \|_{\dot X})^2.$$
But from the previous discussion we already have bounds of the form $\| u_\wb \|_{\dot X} \lesssim 1$ when $\tau$ is small enough.  Thus if we subdivide $[T,T+\tau]$ into small intervals and choose $\eps'$ sufficiently small
as before, we may obtain
$$ \sup_{t \in [T,T+\tau]} \| w \|_{\dot H^1} \lesssim \eps^{10} R^{-10}$$ 
as desired.
\end{proof}

We now apply \eqref{virial-modified} to $u_{\bound,T}$ to obtain
\begin{equation}\label{bound-virial}
\begin{split}
 |\frac{1}{\tau} \int_T^{T+\tau} \int_{B(0,R)} 4 \langle x \rangle^{-\delta} |\nabla u_{\bound,T}|^2 
&- 3 \langle x \rangle^{-\delta} |u_{\bound,T}|^4\ dx dt|\\
&\lesssim \frac{1}{\tau} \sup_{t = T, T+\tau} \int_{B(0,R)} \langle x \rangle^{1-\delta} |u_{\bound,T}(t)| |\nabla u_{\bound,T}(t)|\ dx\\
&+ (\delta + R^{-\delta}) \langle \| u_{\bound,T}\|_{L^\infty_t H^1_x([T, T + \tau] \times \R^3)} \rangle^4.
\end{split}
\end{equation}
We now estimate the various terms in this expression.

First of all, from \eqref{bound-h1} we have
$$ \|u_{\bound,T}(T)\|_{H^1} = \|u_\bound(T)\|_{H^1} \lesssim 1$$
and in particular by mass conservation
$$ \| u_{\bound,T} \|_{L^\infty_t L^2_x([T,T+\tau] \times \R^3)} \lesssim 1.$$
Also from \eqref{bound-close} we have
$$ \| u_{\bound,T} \|_{L^\infty_t \dot H^1_x([T,T+\tau] \times \R^3)} \lesssim 1.$$
Thus we have
$$ \langle \| u_{\bound,T}\|_{L^\infty_t H^1_x([T, T + \tau] \times \R^3)} \rangle^4 \lesssim 1$$
and thus (if we choose $\delta$ small enough depending on $\eps$, and $R$ large enough depending on $\delta$)
$$ (\delta + R^{-\delta}) \langle \| u_{\bound,T}\|_{L^\infty_t H^1_x([T, T + \tau] \times \R^3)} \rangle^4
\lesssim \eps.$$
We now consider the first term in the right-hand side of \eqref{bound-virial}.  Crudely bounding $\langle x \rangle^{1-\delta}$ by $R$ and using the above $H^1$ bounds on $u_{\bound,T}$ we see that this term is bounded by $O(R/\tau)$.  Thus
if we choose $\tau$ large enough depending on $R$ then this term is also $O(\eps)$.  We have thus obtained
\begin{equation}\label{pohozaev-almost}
 |\frac{1}{\tau} \int_T^{T+\tau} \int_{B(0,R)} 4 \langle x \rangle^{-\delta} |\nabla u_{\bound,T}|^2 
- 3 \langle x \rangle^{-\delta} |u_{\bound,T}|^4\ dx dt|
\lesssim \eps.
\end{equation}
We now divide the spatial region of integration $B(0,R)$ on the left into the part $B(0,R_0)$ near the origin,
and the part $B(0,R) \backslash B(0,R_0)$ away from the origin.  First consider the part away from the origin.
From the bounds in \eqref{symbol} (or from Proposition \ref{exterior-decay}), applied for instance with $\delta = 0.1$, 
we have that
\begin{equation}\label{bound-bound}
 |\int_{\R^3 \backslash B(0,R_0)} |\nabla u_\bound|^2 + |u_\bound|^4\ dx| \lesssim R_0^{-2+0.2}
\end{equation}
if $T$ is sufficiently large.  Using \eqref{bound-close} and Sobolev and H\"older, we thus have
$$ |\int_{B(0,R) \backslash B(0,R_0)} |\nabla u_{\bound,T}|^2 + |u_{\bound,T}|^4\ dx| \lesssim R_0^{-2+0.2}$$
if $T$ is sufficiently large.  In particular, if $R_0$ is chosen sufficiently large depending on $\eps$ we have
\begin{equation}\label{bound-bound-T}
 \int_{B(0,R) \backslash B(0,R_0)} |\nabla u_{\bound,T}|^2 + |u_{\bound,T}|^4\ dx \lesssim \eps.
\end{equation}
On the other hand, inside the ball $B(0,R_0)$ we have $\langle x \rangle^{-\delta} = 1 + O(\delta R_0)$ (for instance),
and so again by the $H^1$ bounds on $u_{\bound,T}$ we have
$$\int_{B(0,R_0)} 4 \langle x \rangle^{-\delta} |\nabla u_{\bound,T}|^2 
- 3 \langle x \rangle^{-\delta} |u_{\bound,T}|^4\ dx =
\int_{B(0,R_0)} 4 |\nabla u_{\bound,T}|^2 - 3 |u_{\bound,T}|^4\ dx + O(\delta R_0).$$
Again, if we choose $\delta$ sufficiently small depending on $R_0$ then $O(\delta R_0)$ can be made less than $\eps$.
Combining this with \eqref{bound-bound}, \eqref{bound-bound-T} we thus have
$$ \int_{B(0,R)} 4 \langle x \rangle^{-\delta} |\nabla u_{\bound,T}|^2 
- 3 \langle x \rangle^{-\delta} |u_{\bound,T}|^4\ dx = \int_{\R^3} 4 |\nabla u_\bound|^2 - 3 |u_\bound|^4\ dx
+ O(\eps).$$
Combining this with \eqref{pohozaev-almost} we obtain \eqref{pohozaev-gotit} as desired.
\endprf

{\bf Remark.}  The above argument can be refined slightly to show that
$$
\lim \sup_{T \to +\infty} |\int_T^{T+\tau} \int_{\R^3} 4 |\nabla u_\bound|^2 - 3 |u_\bound|^4\ dx dt|
\lesssim \langle \tau\rangle^\sigma$$
for any $\sigma > 0$, with the implicit constant depending on $\sigma$.  This can be done by setting $R_0$, $R$ to be large
powers of $\tau$, and $\delta$ to be a small power of $\tau$, and using \eqref{symbol} to improve the crude bound of $O(R/\tau)$ obtained above for the first term of \eqref{bound-virial} to $O(R^{0+}/\tau)$.  We omit the details.

{\bf Remark.}  By choosing different multipliers $a$ (not just perturbations of $|x|^2$), it is likely one can
prove other asymptotic identities of this type.  For instance, it seems likely that one has the asymptotic Morawetz
identity
$$
\lim_{\tau \to +\infty}
\lim \sup_{T \to +\infty} |\frac{1}{\tau} \int_T^{T+\tau} (4\pi |u_\bound(t,0)|^2 - \int_{\R^3} \frac{|u_\bound|^4}{|x|}\ dx)| = 0,$$
which (formally at least) is associated to the multiplier $a = |x|$, although the singularity at 0 has to be
treated here with some care.  More generally, it seems that the momentum flux \eqref{flux}, when applied to
$u_\bound$, should converge asymptotically to zero in some weak sense when averaged over increasingly large 
intervals of time.  If one could get stronger control on the decay of the momentum flux, e.g. if it converged
in a suitably strong sense to zero without any averaging in time, this would be substantial progress towards the 
soliton resolution conjecture (solitions have zero momentum current, and are probably the only spherically symmetric solutions to \eqref{nls} with this property).  However, to do this it seems one would first have to discover
a technique that would eliminate the possibility of breather solutions (smooth localized periodic solutions to
\eqref{nls} which are not of the form $u(t,x) = Q(x) e^{iwt}$).  We do not know of a way to rule out such solutions;
the monotonicity formulae arising from virial-type identities do not appear appropriate for this task.  Note that the completely integrable 1D equation supports some breather solutions (see \cite{zakharov}), but this may be a phenomenon
caused by complete integrability and we do not know if there is a similar phenomenon in higher dimensions, say in the spherically symmetric case.


\begin{thebibliography}{10}

\bibitem{bc}
H. Berestycki, T. Cazenave, \emph{Instabilit\'e des \'etats stationnaires dans les \'equations de Schr\"odinger et de Klein-Gordon non lin\'eaires}, C. R. Acad. Sc. Paris, t. \textbf{293} (1981), 489--492. 

\bibitem{bgk}
H. Berestycki, T. Gallou\"et, O. Kavian, \emph{\'Equations de champs scalaires euclidiens non lin'eaires dans le plan.},
C. R. Acad. Sci. Paris S\'er. I Math. \textbf{297} (1983), no. 5, 307--310.

\bibitem{blions}
H. Berestycki, P.L. Lions, \emph{Existence d'ondes solitaires dans des probl\`emes nonlin\'eaires du type Klein-Gordon}, C. R. Acad. Sci. Paris S\'er. A-B \textbf{288} (1979), no. 7, A395--A398.

\bibitem{borg:scatter}
J. Bourgain, \emph{Scattering in the energy space and below for 3D NLS}, J. Anal. Math. \textbf{75} (1998), 267-297. 

\bibitem{borg:growth}
J. Bourgain, \emph{On the growth in time of higher order Sobolev norms of smooth solutions of Hamiltonian PDE}, IMRN \textbf{6} (1996), 277-304. 

\bibitem{borg:book}
J. Bourgain, \emph{New global well-posedness results for non-linear Schr\"odinger equations}, AMS Publications, 1999.

\bibitem{bp}
V. Buslaev, G. Perelman, \emph{Scattering for the nonlinear Schrodinger equations: states close to a solitary wave}, St. Petersburg Math J. \textbf{4} (1993), 1111--1142. 

\bibitem{bp2}
V. Buslaev, G. Perelman, \emph{On the stability of solitary waves for nonlinear Schrodinger equations.  Nonlinear evolution equations} 75-98, Amer. Math. Soc. Transl. Ser. 2, 164, Amer. Math. Soc. Providence, RI 1995. 

\bibitem{bs}
V. Buslaev, C. Sulem, \emph{On asymptotic stability of solitary waves for nonlinear Schrodinger equations}, Ann. Inst. H. Poincare Anal. Nonlineaire 20 (2003) 3, 419-475. 

\bibitem{caz}
T. Cazenave, An introduction to nonlinear Schrodinger equations,  Textos de Metodes Matematicos \textbf{22} (Rio de Janeiro), 1989. 

\bibitem{caz-lions}
T. Cazenave, P. Lions, \emph{Orbital stability of standing waves for some nonlinear Schrodinger equations}, Comm. Math. Phys. \textbf{68} (1979), 209--243. 

\bibitem{cwI}
T. Cazenave, F.B. Weissler, \emph{Critical nonlinear Schr\"odinger
Equation}, Non. Anal. TMA, \textbf{14} (1990), 807--836.

\bibitem{christ}
M. Christ, Lectures on singular integral operators. CBMS Regional Conference Series in Mathematics, 77. 
 American Mathematical Society, Providence, RI, 1990. 

\bibitem{coff}
C.V. Coffman, \emph{Uniqueness of the ground state solution for $\Delta u - u + u^3 = 0$ and a variational characterization of other solutions}, Arch. Rat. Mech. Anal. \textbf{46} (1972), 81--95.

\bibitem{ckstt:scatter}
J. Colliander, M. Keel, G. Staffilani, H. Takaoka, T. Tao, \emph{Scattering for the 3D cubic NLS below the energy norm}, to appear, Comm. Pure Appl. Math.

\bibitem{ckstt:7}
J. Colliander, M. Keel, G. Staffilani, H. Takaoka, T. Tao, \emph{Polynomial growth and orbital instability bounds for the 1D cubic NLS below the energy norm}, Discrete Cont. Dynam. Systems. \textbf{2} (2003), 33--50.

\bibitem{ckstt:8}
J. Colliander, M. Keel, G. Staffilani, H. Takaoka, T. Tao, \emph{ Polynomial growth and orbital instability bounds for $L^2$-subcritical NLS below the energy norm }, to appear, Comm. Pure Appl. Anal.

\bibitem{cuccagna}
S. Cuccagna, \emph{Stabilization of solutions to nonlinear Schrodinger equations}, CPAM \textbf{54} (2001), 1110--1145. 

\bibitem{cuccagna-2}
S. Cuccagna, \emph{On asymptotic stability of ground states of NLS}, preprint. 

\bibitem{gv:scatter}
J. Ginibre, G. Velo, \emph{Scattering theory in the energy space for a class of nonlinear Schr\"odinger equations}, J. Math. Pure. Appl. \textbf{64} (1985), 363--401.

\bibitem{glassey}
R.T. Glassey, \emph{On the blowing up of solutions to the Cauchy problem for nonlinear Schrodinger operators}, J. Math. Phys. \textbf{8} (1977), 1794--1797. 

\bibitem{gss}
M. Grillakis, J. Shatah, W. Strauss, \emph{Stability theory of solitary waves in the presence of symmetry I}, J. Funct. Anal. \textbf{74} (1987), 160--197. 

\bibitem{gss-2}
M. Grillakis, J. Shatah, W. Strauss, \emph{Stability theory of solitary waves in the presence of symmetry II}, J. Funct. Anal.94 (1990), 308--348. 

\bibitem{tao:keel}
M. Keel, T. Tao, \emph{Endpoint Strichartz Estimates}, Amer. Math. J. 120 (1998), 955--980.

\bibitem{kpv}
C. Kenig, G. Ponce, L. Vega, \emph{Well-posedness and scattering results for the generalized Korteweg-de Vries equations via the contraction principle}, Comm. Pure Appl. Math. \textbf{46} (1993), 527--620.

\bibitem{mmt}
Y. Martel, F. Merle, T-P. Tsai, \emph{Stability and asymptotic stability in the energy space of the sum of N solitons for subcritical gKdV equations}, preprint.

\bibitem{merle1} F. Merle, \emph{Construction of solutions with exactly k blow-up points for the Schrodinger equation with critical non-linearity}, Comm. Math. Phys. \textbf{149} (1992), 205--214. 

\bibitem{merle2} F. Merle, \emph{Asymptotics for $L^2$ minimal blowup solutions of critical nonlinear Schrodinger equation}, Ann. Inst. Henri Poincare 13 (1996), 553--565. 

\bibitem{nak:scatter}
K. Nakanishi, \emph{Energy scattering for non-linear Klein-Gordon and Schrodinger equations in spatial dimensions 1 and 2}, JFA \textbf{169} (1999), 201--225. 

\bibitem{ogawa}
T. Ogawa, Y. Tsutsumi, \emph{Blow-up of $H^1$ solution for the nonlinear Schr\"odinger equation}, J. Differential Equations \textbf{92} (1991), no. 2, 317-330. 

\bibitem{ozawa} 
T. Ozawa, \emph{Long range scattering for nonlinear Schrodinger equations in one space dimension}, CMP \textbf{139} (1991), 479--493.

\bibitem{perelman}
G. Perelman, \emph{Some results on the scattering of weakly interacting solitons for nonlinear Schr\"odinger equations}, in ``Spectral theory, microlocal analysis, and singular manifolds'', Akad. Verlag. (1997), 78--137.

\bibitem{rauch}
J. Rauch, \emph{Local decay of scattering solutions to Schrodinger's equation}, CMP \textbf{61} (1978), 149--168. 

\bibitem{schlag}
I. Rodnianski, W. Schlag, A. Soffer, \emph{Asymptotic stability of $N$-soliton states of NLS}, preprint.

\bibitem{igor}
I. Rodnianski, T. Tao, \emph{Long-time decay estimates for the Schr\"odinger equation on compact perturbations of Euclidean space},  preprint.

\bibitem{segur}
H. Segur, M. Ablowitz, \emph{Asymptotic solutions and conservation laws for the nonlinear Schr\"odinger equation I}. J. Math. Phys. \textbf{17} (1976), 710--713. 

\bibitem{shatah}
J. Shatah, W. Strauss, \emph{Instability of nonlinear bound states}, Comm. Math. Phys. \textbf{100} (1985), 173--190.

\bibitem{sjolin}
P. Sjolin, \emph{Regularity of solutions to the Schrodinger equation}, Duke Math. J. \textbf{55} (1987), 699-715. 

\bibitem{staff:growth}
G. Staffilani, \emph{On the growth of high Sobolev norms of solutions for KdV and Schrodinger equations}, Duke Math J. \textbf{86} (1997), 109-142. 

\bibitem{strauss:soliton}
W. Strauss, \emph{Existence of solitary waves in higher dimensions}, Comm. Math. Phys. \textbf{55} (1977), 149--162.

\bibitem{strauss}
W. Strauss, Nonlinear wave equations, Regional Conf. Series in Math., 1989. 

\bibitem{stein:large}
E.~M. Stein, \emph{Harmonic Analysis}, Princeton University Press, 1993.

\bibitem{heron}
R. Strichartz, \emph{Asymptotic behavior of waves}, J. Funct. Anal. \textbf{40} (1981), 341-357.

\bibitem{sulem-break}
C. Sulem, P. Sulem, \emph{Focusing nonlinear Schr\"odinger equation and wave packet collapse}, Proceedings of the Second World Congress of Nonlinear Analysts, Part 2 (Athens, 1996). Nonlinear Anal. \textbf{30} (1997), no. 2, 833--844.

\bibitem{sulem}
C. Sulem, P. Sulem, The nonlinear Schrodinger equation: Self-Focusing and Wave Collapse, Applied Mathematical Sciences 139, Springer-Verlag, New York.

\bibitem{tsai-yau}
T.P. Tsai, H.T. Yau, \emph{Asymptotic dynamics of nonlinear Schrodinger equations: resonance dominated and dispersion dominated solutions}, CPAM \textbf{55} (2002), 153--216. 

\bibitem{tsai}
T.P. Tsai, \emph{Asymptotic dynamics of nonlinear Schrodinger equations with many bound states}, JDE \textbf{192} (2003), 225--282. 

\bibitem{vega}
L. Vega, \emph{Schrodinger equations: pointwise convergence to the initial data}, Proc. Amer. Math. Soc. \textbf{
102} (1988), 874-878. 

\bibitem{vilela}
M. Vilela, \emph{Regularity of solutions to the free Schrodinger equation with radial initial data}, Ill. J. Math. \textbf{45} (2001), 361--370. 

\bibitem{wein1}
M. Weinstein, \emph{Nonlinear Schrodinger equations and sharp interpolation estimates}, Comm. Math. Phys. \textbf{87} (1983), 567--576. 

\bibitem{wein2}
M. Weinstein, \emph{Modulational stability of ground states of nonlinear Schrodinger equations}, SIAM J. Math. Anal. \textbf{16} (1985), 472--491. 

\bibitem{zakharov}
V.E. Zakharov, A.B. Shabat, \emph{Exact theory of two-dimensional self-focusing and one-dimensional
self-modulation of waves in nonlinear media}, Soviet Physics JETP v34 no. 1, (1972) 62--69.

\end{thebibliography}
\end{document}